\documentclass{article}
\usepackage[utf8]{inputenc}
\usepackage{amsmath}
\usepackage{amsfonts}
\usepackage{amsthm}
\usepackage{tikz}
\usepackage{xcolor}
\usepackage{newtxtext,newtxmath}
\usepackage{caption}
\usepackage{subcaption}
\usepackage{hyperref}
\usepackage{authblk}
\usepackage{todonotes}
\usepackage{soul}

\theoremstyle{definition}

\theoremstyle{remark}
\newtheorem*{remark}{Remark}

\newcommand{\reviewerA}[1]{{\color{black}#1}}
\newcommand{\reviewerB}[1]{{\color{black}#1}}

\newcommand{\both}[1]{{\color{black}}{#1}}
    
\def\tolopt{\text{tol}_{\text{opt}}}
\def\maxit{\text{it}_{\text{max}}}
\def\nmax{N_{\text{max}}}

\title{An optimisation--based domain--decomposition reduced order model for parameter-dependent non--stationary fluid dynamics problems}

\date{}
\author[a]{Ivan Prusak\footnote{\href{mailto:iprusak@sissa.it}{iprusak@sissa.it}}}
\author[a]{Davide Torlo\footnote{\href{mailto:dtorlo@sissa.it}{dtorlo@sissa.it}}}
\author[b]{Monica Nonino\footnote{\href{mailto:monica.nonino@univie.ac.at}{monica.nonino@univie.ac.at}}}
\author[a]{\\Gianluigi Rozza\footnote{\href{mailto:grozza@sissa.it}{grozza@sissa.it}}} 

\affil[a]{Mathematics Area, mathLab, SISSA, 34136 Trieste, Italy}
\affil[b]{Fakult\"at f\"ur Mathematik,
Universit\"at Wien, 1090 Wien, Austria}

\graphicspath{../}

\begin{document}

\maketitle

\begin{abstract}
\label{abstact}
In this work, we address parametric non--stationary fluid dynamics problems within a model order reduction setting based on domain decomposition. 
Starting from the optimisation--based domain decomposition approach, we derive an optimal control problem, for which we present a convergence analysis in the case of non--stationary incompressible Navier--Stokes equations. 
\reviewerB{We discretize the problem with the finite element method and we compare different model order reduction techniques: POD--Galerkin and a non--intrusive neural network procedures.
We show that the classical POD--Galerkin is more robust and accurate also in transient areas, while the neural network can obtain simulations very quickly though being less precise in the presence of discontinuities in time or parameter domain.}
We test the proposed methodologies on two fluid dynamics benchmarks \reviewerB{with physical parameters and time dependency}: the non--stationary backward--facing step and lid--driven cavity flow. 
\end{abstract}


\section{Introduction}
\label{introduction}

With the increase in the potential of high--performance computing in the last years, there is a significant demand for numerical methods and approximation techniques that can perform real--time simulations of Partial Differential Equation (PDEs). The applications vary from naval, aeronautical and biomedical engineering.  There exist many techniques to achieve such a goal, including reduced--order modelling~\cite{Rozza_book} and domain--decomposition (DD) methods~\cite{QuarteroniValiDD}. 

The DD methodology is a highly efficient tool in the framework of PDEs. Any DD algorithm is constructed by an effective splitting of the domain of interest into different subdomains (overlapping or not), and the original problem is then restricted to each of these subdomains with some coupling conditions on the intersections of the subdomains. The coupling conditions may be very different, they depend on the physical meaning of the problem at hand, and they are required to provide a certain degree of continuity among these subdomains (see, for example,  \cite{QuarteroniValiNumerics, QuarteroniValiDD}).  These methods are extremely important for multi--physics problems when efficient subcomponent numerical codes are already available, or when we do not have direct access to the numerical algorithms for some parts of the systems; see, for instance, \cite{ErvinJenkinsLee2014, GosseletChiaruttiniReyFeyel2012, HoangLee2021, KUBERRY2013594, Kuberry2015, lagnese2004domain}. 

Model--order reduction methods are another set of methods that are extremely useful when dealing with real--time simulations or multi--query tasks. These methods are successfully employed in the settings of non--stationary and/or parameter--dependent PDEs. Reduced--order models (ROMs) are extremely effective due to the splitting of the computational effort into two stages: the offline stage, which contains the most expensive part of the computations, and the online stage, which allows performing fast computational queries using structures that are pre--computed in the offline stage; for more details, see \cite{Rozza_book}. ROMs have been successfully applied in different fields such as fluid dynamics \cite{ALI20202399,carere2021weighted,crisovan2019model, DeparisRozza2009, LassilaManzoniQuarteroniRozza2014, Rozza2009ReducedBM,stabile2019reduced,stabile2018finite,strazzullo2020pod,strazzullo2022consistency,tezzele2020enhancing,torlo2021model}, structural mechanics \cite{ballarin2016pod,ballarin2017reduced, Haasdonk, RozzaHuynhPatera2007,torlo2018stabilized,venturi2019weighted} and fluid--structure interaction problems \cite{AstorinoChoulyFernandez2010, ballarin2016pod, Nonino2021, Nonino2020}. Lately, there have also been great advances in reduced--order modelling for optimal--control problems, \cite{pichi2022driving, strazzullo2018model, strazzullo2022pod}.

\reviewerB{The novelty of this paper lies in the study of an optimization-based approach for a domain-decomposition model addressing unsteady parametric fluid dynamics problems. We employ both intrusive and non-intrusive model-order reduction (MOR) techniques to achieve this, providing a comprehensive comparison between the two methodologies. Utilizing domain-decomposition methods, particularly an optimization approach ensuring coupling of interface conditions between subdomains as presented in \cite{GunzburgerLee2000, GUNZBURGER2000177, prusak2022optimisationbased}, we develop various reduced-order models (ROMs). These include classical POD-Galerkin projection-based intrusive methods \cite{Rozza_book, stabile2018finite} and neural network (NN) based non-intrusive ROMs \cite{hesthaven2018non, romor2023non, papapicco2022neural, pichi2021artificial, siena2023fast}. Moreover, we perform a convergence analysis for the nonstationary fluid dynamics problems based on the \textit{a priori} estimates by extending the results presented in \cite{GunzburgerLee2000} for the stationary case.}


This work is structured as follows. In Section~\ref{problem_formulation}, we introduce the monolithic fluid dynamics problem and its time--discretisation scheme with the further derivation of the optimisation--based domain--decomposition formulation at each time step. In Section~\ref{analysis}, we derive \textit{a priori} estimates for the solutions to Navier--Stokes equations which are then used to prove the existence and uniqueness of the minimiser to the optimal--control problem derived in the previous section.  \both{Furthermore, in Section~\ref{optimality_system} we provide a closed--form expression for the gradient of the objective functional which allows for the use of gradient--based optimisation algorithms to decouple the subdomain solves. } Section~\ref{high_fidelity} contains the Finite Element discretisation of the problem of interest and the corresponding finite--dimensional high--fidelity optimisation problem. Section~\ref{ROM} deals with two ROM techniques: an intrusive Galerkin projection and a neural network (NN) algorithm, both based on a Proper Orthogonal Decomposition (POD) methodology. In Section~\ref{results}, we show some numerical results for two toy problems: the backward--facing step and the lid--driven cavity flow. Conclusions will follow in Section~\ref{conclusions}.

\section{Problem formulation}
\label{problem_formulation}
In this section, starting with a monolithic formulation of the time--dependent incompressible Navier--Stokes equations, we first introduce a time discretisation on the continuous level employing the implicit Euler time--stepping scheme. Then, we will describe a two--domain optimisation-based domain--decomposition formulation at each time step in variational form in the end. Here and in the next few sections, the analysis is valid for any value of the physical parameter, so for the sake of simplicity, we postpone mentioning the parameter dependence of the problem until Section~\ref{ROM}.

\subsection{Monolithic formulation}
\label{monolithic_formulation}

\begin{figure}\centering
\begin{subfigure}{0.49\textwidth}
        \centering
    \begin{tikzpicture}
 \draw[very thick] (0,0) ellipse (1.5cm and 2cm);
\node at (0,0) {$\Omega$};
\draw (0,1.9) -- (0,2.1);
\draw (0,-1.9) -- (0,-2.1);
\node[anchor=east] at (-1.6,0) {$\Gamma_{D}$};
\node[anchor=west] at (1.6,0) {$\Gamma_{N}$};
\end{tikzpicture}
    \caption{Physical domain}
    \label{fig:mono_domain}
\end{subfigure}\hfill
\begin{subfigure}{0.49\textwidth}
\centering
    \begin{tikzpicture}
    
\draw[red, very thick]  (1.07121,1.4) parabola (-1.45236,-0.5);
\node[red, very thick, anchor=east] at (-0., 0.3) {\textbf{$\Gamma_0$}};
\node at (-0.5,1.2) {$\Omega_1$};
\node at (0.1,-0.5) {$\Omega_2$};
\draw (0,1.9) -- (0,2.1);
\draw (0,-1.9) -- (0,-2.1);
\draw (1.17,1.47) -- (0.97,1.33);
\draw (-1.55236,-0.53) -- (-1.35236,-0.47);

\node[anchor=east, color=magenta] at (-1.4,0.7) {$\Gamma_{D,1}$};
\node[anchor=west, color=teal] at (1.5,0) {$\Gamma_{N,2}$};
\node[anchor=west, color=olive] at (0.68,1.85) {$\Gamma_{N,1}$};
\node[anchor=east, color=violet] at (-1,-1.5) {$\Gamma_{D,2}$};

\draw[domain=0:1.07121, smooth, very thick, variable=\x, samples=51, olive] plot ({\x}, {2*sqrt(1.0-\x*\x/2.25)});

\draw[domain=-2:1.4, smooth, very thick, variable=\y, samples=501, teal] plot ({sqrt(2.25*(1-\y*\y/4))}, {\y});

\draw[domain=-1.45236:0, smooth, very thick, variable=\x, samples=51, violet] plot ({\x}, {-2*sqrt(1.0-\x*\x/2.25)});

\draw[domain=-0.5:2, smooth, very thick, variable=\y, samples=501, magenta] plot ( {-sqrt(2.25*(1-\y*\y/4))},{\y});

\end{tikzpicture}
    \caption{Domain Decomposition}
    \label{fig:dd_domain}
    \end{subfigure}
    \caption{Domain and boundaries}
\end{figure}
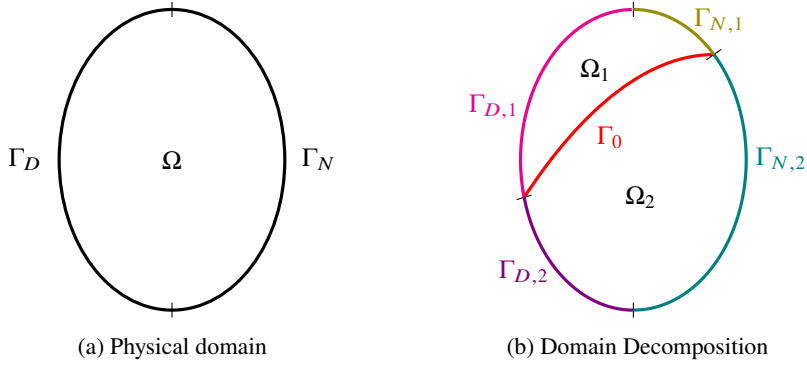

Let $\Omega$ be a physical domain of interest: we assume $\Omega$ to be an open subset of $\mathbb{R}^2$ and $\Gamma$ to be the boundary of $\Omega$. We also consider a finite time interval $[0, T]$ with $T >0$.
Let $f: \Omega \times [0,T] \rightarrow \mathbb{R}^2$ be the forcing term, $\nu$ the kinematic viscosity, $u_{D}$ a given Dirichlet datum and $u_0$ a given initial condition.  The problem reads as follows: find the velocity field $u: \Omega\times [0,T] \rightarrow \mathbb{R}^2 $ and the pressure $p: \Omega\times [0,T] \rightarrow \mathbb{R}$ s.t.

\begin{subequations}   \label{eq:mono}
\begin{eqnarray}
    \frac{\partial u}{\partial t} -\nu \Delta u + \left( u \cdot \nabla \right) u + \nabla p = f  & \text{in} & \Omega\times (0,T], \label{eq:mono1}\\
    -\text{div} u = 0& \text{in} & \Omega\times (0,T], \label{eq:mono2}  \\
    u = u_{D}  & \text{on} & \Gamma_{D}\times [0,T], \label{eq:mono3} \\
   \nu \frac{\partial u}{\partial \textbf{n}} - p \textbf{n} = 0 & \text{on} & \Gamma_{N}\times [0,T] \label{eq:mono4},
   \\ u(t=0) = u_0 & \text{in} & \Omega \label{eq:mono5},
\end{eqnarray}
\end{subequations}
where $\Gamma_{D}$ and $\Gamma_{N}$ are  disjoint subsets of $\Gamma$ (as it is shown in Figure \ref{fig:mono_domain}) and $\textbf{n}$ is an outward unit normal vector to $\Gamma_{N}$.

\subsection{Time discretisation}

We will start with time discretisation of problem \eqref{eq:mono}. Let $\Delta t >0$, we assume the following time interval partition: $0 = t_0 < t_1 < .... < t_M =T$, where $t_n = n\Delta t$ for $n = 0, ..., M$. We employ the implicit Euler scheme for the incompressible Navier--Stokes equation which reads as follows: for $n \geq 1$ find $u^n: \Omega \rightarrow \mathbb{R}^2$, $p^n: \Omega \rightarrow \mathbb{R}$  s.t. 
\begin{subequations}   \label{eq:TD_mono_equation}
\begin{eqnarray}
    \frac{u^{n} - u^{n-1}}{\Delta t} -\nu \Delta u^{n} + \left( u^{n} \cdot \nabla \right) u^{n} + \nabla p^{n} = f^{n}  & \text{in} & \Omega, \label{eq:TD_mono1}\\
    -\text{div} u^{n} = 0& \text{in} & \Omega, \label{eq:TD_mono2}  \\
    u^{n} = u_{D}^{n}  & \text{on} & \Gamma_{D}, \label{eq:TD_mono3} \\
   \nu \frac{\partial u^{n}}{\partial \textbf{n}} - p^{n} \textbf{n} = 0 & \text{on} & \Gamma_{N} \label{eq:TD_mono4},
\end{eqnarray}
and for $n=0$
\begin{equation}
    u^0 = u_0 \quad \text{in} \quad \Omega. \label{eq: TD_mono0}
\end{equation}
\end{subequations}
Here we adopted the following notations: $f^n(\cdot) = f(\cdot, t_n)$ and $u^n_D(\cdot) = u_D(\cdot, t_n)$.

\subsection{Domain Decomposition (DD) formulation} 
\label{DD}
\reviewerB{As mentioned in the introduction, we resort to the optimisation--based DD framework in \cite{prusak2022optimisationbased}.} Let $\Omega_i, \ i=1,2$, be open subsets of $\Omega$, such that  $\overline{\Omega} = \overline{\Omega_1 \cup \Omega_2}$, \  $\Omega_1 \cap \Omega_2 = \emptyset$. Denote $\Gamma_i := \partial \Omega_i \cap \Gamma, \ i=1,2,$ and $\Gamma_0 := \overline{\Omega_1} \cap \overline{\Omega_2}$. In the same way, we define the corresponding boundary subsets $\Gamma_{i,D}$ and $\Gamma_{i, N}$, $i=1,2$, see Figure~\ref{fig:dd_domain}. 
\reviewerB{Having at hand the time discretisation~\eqref{eq:TD_mono_equation}, we can cast the DD problem in the optimal control setting with the weak formulation of the DD state equations. For this purpose, we define the following variational spaces  
\begin{itemize}
    \item $V_i := \left\{ u \in H^1(\Omega_i; \mathbb{R}^2)\right\}$,
    \item $V_{i,0} := \left\{ u \in H^1(\Omega_i; \mathbb{R}^2): \ u|_{\Gamma_{i, D}  } =0\right\}$,
    \item $Q_i := \left\{ p \in L^2(\Omega_i; \mathbb{R}) \right\}$,
\end{itemize} 
together with the following bilinear and trilinear forms
\begin{itemize}
    \item  $m_i: V_i \times V_{i,0} \rightarrow \mathbb{R}, \quad m_i(u_i,v_i) :=  (u_i,  v_i)_{\Omega_i}$,
    \item  $a_i: V_i \times V_{i,0} \rightarrow \mathbb{R}, \quad a_i(u_i,v_i) := \nu (\nabla u_i, \nabla v_i)_{\Omega_i} $,
     \item  $b_i: V_i \times Q_{i} \rightarrow \mathbb{R}, \quad b_i(v_i,q_i) := - (\text{div} v_i, q_i)_{\Omega_i}$,
     \item  $c_i: V_i \times V_i \times V_{i,0} \rightarrow \mathbb{R}, \quad c_i(u_i,w_i, v_i) := \left( (u_i \cdot \nabla)w_i, v_i\right)_{\Omega_i}   $,
\end{itemize}
for $i=1,2$. In the definitions above, $(\cdot, \cdot)_{\omega}$ indicates the $L^2(\omega)$ inner product. The spaces $V_i$ are endowed with the $H^1(\Omega_i)$--norm for $i=1,2$, the spaces $V_{i,0}$ with the $H_0^1(\Omega_i)$--norm and the spaces $Q_i$ with the $L^2(\Omega_i)$--norm for $i=1,2$. 
}

\noindent\reviewerB{The DD formulation then reads as follows. For $n \geq 1$ solve the following optimal control  problem: 
\begin{align}
&\text{\reviewerB{ \textit{minimise over } $g \in L^2(\Gamma_0)$ \textit{ the functional}} } \nonumber\\
  & \label{eq:functional}  \mathcal J_\gamma(u_1^n, u_2^n; g) := \frac{1}{2} \int_{\Gamma_0} \left| u_1^n - u_2^n\right|^2 d\Gamma + \frac{\gamma}{2}\int_{\Gamma_0} \left|g\right|^2 d\Gamma, \\
  &\nonumber \textit{subject to the variational problem:}\\
  &\begin{cases}
     \label{eq:state}
        \frac{1}{\Delta t} m_i \left( u_i^n - u_i^{n-1}, v_i\right)+ a_i(u_i^n, v_i) + c_i(u_i^n, u_i^n, v_i)  + b_i(v_i, p_i^n)    &\\ 
 \qquad\qquad\qquad\qquad  =   (f_i^n, v_i)_{\Omega_i}
     + \left( (-1)^{i+1} g, v_i \right)_{\Gamma_0} \quad \quad&\forall v_i \in V_{i,0},   \\
    b_i(u_i^n, q_i)   = 0 \quad \quad  &\forall q _i \in Q_{i},  \\ 
     u_i^n    =   u_{i, D}^n  \quad \quad & \text{on}  \ \Gamma_{i, D}.
  \end{cases}
\end{align}   
}

\reviewerB{The presence of the regularisation term $\frac{\gamma}{2}\int_{\Gamma_0} \left|g\right|^2 d\Gamma$ with $\gamma>0$ in the definition of the objective functional~\eqref{eq:functional} is essential for the well--posedness analysis of the optimal control problem~\eqref{eq:functional}--~\eqref{eq:state} carried out in the next section. Nevertheless, our goal is to find a DD solution which is as close as possible to the monolithic one; for this reason we will provide the analysis when the regularisation parameter $\gamma \rightarrow 0$.     }

\section{Analysis of the optimal control problem}\label{analysis}
In this section, we will give an overview of the existence of local minima of the optimal--control problem described above. 
It will rely on the a priori estimates for the solutions to the Navier--Stokes equations.
\reviewerB{As is evidenced, for instance, in \cite{Richter2017}, due to the presence of the Neumann boundary condition, the analysis of the state problem is not possible in the general case}, so we will modify the problem in the framework where the problem is well posed and give some indication about the original problem later in the section. \both{Finally, in the end of the section, we will list the objective functional gradient of the original DD optimal control problem.  }

\subsection{A modified Navier--Stokes problem}

First, without loss of generality, we assume that the Dirichlet data $u_{i, D}$ is homogeneous. Otherwise, it is possible to recast the inhomogeneous problem into a homogeneous one with various techniques, e.g., lifting functions as \both{described in Appendix~\ref{lifting_supremiser}}.
As mentioned in the preface to this section, it is hard to prove the well--posedness of the solution to the Navier--Stokes equation in the form \eqref{eq:state}. The main problem arises from the nonlinear term $c_i(u_i^n, u_i^n,v_i) = \left( (u_i^n \cdot \nabla)u_i^n, v_i \right)_{\Omega_i^n}$. Indeed, by integration by parts and the incompressibility conditions in \eqref{eq:state}, we can see
\reviewerB{
\begin{eqnarray*}
    c_i(u_i^n, w_i, v_i) & & = \int_{\Omega_i} (u_i^n \cdot \nabla) w_i\cdot v_i d\Omega = \int_{\partial \Omega_i} (w_i \cdot v_i) (u_i^n\cdot \textbf{n}_i)d\Gamma 
    \\ & & - \int_{\Omega_i} (w_i \cdot v_i) \,\text{div}\, u_i^n d\Omega - \int_{\Omega_i} (u_i^n \cdot \nabla )v_i \cdot w_i d\Omega 
    \\ & & = \int_{\partial \Omega_i} (w_i \cdot v_i) (u_i^n\cdot \textbf{n}_i)d\Gamma  - \int_{\Omega_i} (u_i^n \cdot \nabla )v_i \cdot w_i d\Omega 
    \\ & & =  \int_{\Gamma_{i,N} \cup \Gamma_0} (w_i \cdot v_i) (u_i^n\cdot \textbf{n}_i)d\Gamma  - c_i(u_i^n, v_i, w_i),
\end{eqnarray*}
}
which leads to the following expression 
\begin{equation}
    c_i(u_i^n, u_i^n, u_i^n) = \frac{1}{2}  \int_{\Gamma_{i,N} \cup \Gamma_0} |u_i^n|^2  (u_i^n\cdot \textbf{n}_i)d\Gamma . \label{eq:boundary_term_c}
\end{equation}
\reviewerB{The issue here is that this boundary term is of third order and has an unknown sign since the solution on the interface $\Gamma_0$ is undetermined. This complicates the analysis, which relies on the \textit{a priori} bounds of the weak solutions.
On the other hand, it gives us an idea of how to redefine the problem at hand in order to obtain well--posedness (see, e.g. \cite{M2AN_1996__30_7_815_0, Richter2017}). Indeed,  the trilinear forms $c_i(\cdot, \cdot, \cdot), i=1,2$ can be rewritten as 
\begin{eqnarray*}
    c_i(u_i, w_i, v_i) & = & \frac{1}{2}c_i(u_i, w_i, v_i) + \frac{1}{2}c_i(u_i, w_i, v_i) 
    \\ & = & \frac{1}{2}c_i(u_i, w_i, v_i) - \frac{1}{2}c_i(u_i, v_i, w_i) + \frac{1}{2}\int_{\Gamma_{i,N} \cup \Gamma_0} (w_i \cdot v_i) (u_i\cdot \textbf{n}_i)d\Gamma 
    \\ & = & \tilde c_i(u_i, w_i, v_i) + \frac{1}{2}\int_{\Gamma_{i,N} \cup \Gamma_0} (w_i \cdot v_i) (u_i\cdot \textbf{n}_i)d\Gamma,
\end{eqnarray*}
where the trilinear form $\tilde c (\cdot, \cdot, \cdot)$ is defined as 
\begin{equation}
    \label{eq:trilinear_new} \tilde c_i(u_i, w_i, v_i) = \frac{1}{2} \left( (u_i \cdot \nabla)w_i, v_i\right)_{\Omega_i} - \frac{1}{2} \left( (u_i \cdot \nabla)v_i, w_i\right)_{\Omega_i}
\end{equation}
and it has the following remarkable property
\begin{eqnarray}\label{eq:trilinear_property}
\tilde c_i(u_i, v_i, v_i) = 0 \quad \forall u_i, v_i \in V_{i,0}.
\end{eqnarray}
We thus consider a new variational formulation of the state equations~\eqref{eq:state}: for $n\geq 1$ and $i=1,2$, find $u_i^n \in V_{i,0}$ and $p_i^n \in Q_i$ s.t.
\begin{subequations}\label{eq:state_new}
\begin{align}
   \label{eq:state1_new}   \begin{split}
          \frac{1}{\Delta t} m_i \left( u_i^n - u_i^{n-1}, v_i\right)+ &a_i(u_i^n, v_i) + \tilde c_i(u_i^n, u_i^n, v_i)  + b_i(v_i, p_i^n)    \\ 
 & \quad  =   (f_i^n, v_i)_{\Omega_i}
     + \left( (-1)^{i+1} g, v_i \right)_{\Gamma_0} 
  \quad \quad    \forall v_i \in V_{i,0},         \end{split}&\\
  \label{eq:state2_new} & b_i(u_i^n, q_i)   = 0 \quad \quad  \forall q _i \in Q_{i}, & 
\end{align}
\end{subequations}
which corresponds to imposing the following modified Neumann outlet condition
\begin{eqnarray} \label{eq:outlet_new}
    \nu \frac{\partial u_i^n}{\partial \textbf{n}_i} - p_i^n \textbf{n}_i - \frac{1}{2} (u_i^n \cdot \textbf{n}_i) u_i^n = 0 & \text{on} & \Gamma_{i,N}.
\end{eqnarray} 
The boundary condition~\eqref{eq:outlet_new} no longer arises from the original physical problem~\eqref{eq:mono}, but, on the other hand, the condition~\eqref{eq:trilinear_property} allows us to derive the \textit{a priori} estimates in the next section. 
It is also worth mentioning that in the case of DD optimal control problem~\eqref{eq:functional}--\eqref{eq:state_new} the optimum $g$ is the approximation of the modified normal stress on the interface, i.e., $\nu \frac{\partial u_i^n}{\partial \textbf{n}_i} - p_i^n \textbf{n}_i - \frac{1}{2} (u_i^n \cdot \textbf{n}_i) u_i^n,\, i=1,2,$ whereas in the case of the original problem~\eqref{eq:functional}--\eqref{eq:state} the optimum $g$  is the approximation of the physical normal stress $\nu \frac{\partial u_i^n}{\partial \textbf{n}_i} - p_i^n \textbf{n}_i, i=1,2$. }

\subsection{\textit{A priori} estimates}
\reviewerB{The key element for the analysis of the optimal control problem is the derivation of the \textit{a priori} estimates of the solutions to the state equations~\eqref{eq:state_new}. We will follow the idea introduced in~\cite{GunzburgerLee2000} for the stationary Navier--Stokes equations and extend the results to the non--stationary case.
We first introduce the various well--known properties of the different terms in the weak formulation \eqref{eq:state_new}, see, for instance,~\cite{GunzburgerLee2000, Richter2017}.}
\begin{itemize}
    \item The forms $m_i(\cdot, \cdot)$, $a_i(\cdot, \cdot)$ and $\tilde c_i (\cdot,\cdot,\cdot)$ are continuous: there exist positive constants $C_m$, $C_a$ and $C_c$ such that
    \begin{align}
    \label{eq:m_cont} & |m_i(u_i, v_i)|  \leq C_m || u_i ||_{V_{i,0}} || v_i ||_{V_{i,0}} \quad &\forall u_i, v_i \in V_{i,0},
    \\   & \label{eq:a_cont} |a_i(u_i, v_i)|  \leq  C_a || u_i ||_{V_{i,0}} || v_i ||_{V_{i,0}} \quad &\forall u_i, v_i \in V_{i,0},
        \\\label{eq:c_cont} & |\tilde c_i(u_i, w_i,v_i)| \leq  C_c || u_i ||_{V_{i,0}} || w_i ||_{V_{i,0}} || v_i ||_{V_{i,0}} \quad &\forall u_i,w_i, v_i \in V_{i,0},
    \end{align}
\item the bilinear form  $a_i(\cdot, \cdot)$ is coercive: there exists a positive constant $\alpha >0$ such that
\begin{eqnarray}
    \label{eq:a_coercive} a_i(v_i, v_i) \geq \alpha ||v_i||_{V_{i,0}}^2 \quad \forall v_i \in V_{i,0},
\end{eqnarray}
\item the bilinear form  $b_i(\cdot, \cdot)$ satisfies inf--sup condition: there exists a positive constant $\beta >0$ such that
\begin{eqnarray}
    \label{eq:b_infsup} \sup\limits_{v_i \in V_{i,0}\backslash\{0\}} \frac{b_i(v_i, q_i)}{||v_i||_{V_{i,0}}} \geq \beta ||q_i||_{Q_{i}} \quad \forall   q_i \in Q_i,
\end{eqnarray}
\item the bilinear form $m_i(\cdot, \cdot)$ is non--negative definite, i.e. 
    \begin{eqnarray}
        \label{eq:m_nonneg} m_i(v_i, v_i) = || v_i ||_{L^2(\Omega_i)}^2 \geq 0 \quad \forall v_i \in V_{i,0}.
    \end{eqnarray}
By using the properties \eqref{eq:m_cont}, \eqref{eq:a_coercive}, \eqref{eq:trilinear_property}, \eqref{eq:m_nonneg}, the trace theorem and  equations \eqref{eq:state_new}, we are able to write the following estimate for the solution $u_i^n$ and $p_i^n$ to~\eqref{eq:state_new}
\begin{eqnarray*}
    ||u_i^n||_{V_{i,0}}^2 &  \leq & \frac{1}{\alpha} a_i(u_i^n, u_i^n) \leq\frac{1}{\alpha} a_i(u_i^n, u_i^n) + \frac{1}{\alpha\Delta t} m_i(u_i^n, u_i^n) 
 \\  & = & \frac{1}{\alpha} \left( \frac{1}{\Delta t} m_i(u_i^{n-1}, u_i^n) - \tilde c_i(u_i^n, u_i^n, u_i^n) - b_i(u_i^n, p_i^n) \right.
 \\ & +& \left. (f_i^n, u_i^n)_{\Omega_i} + (-1)^{i+1}(g, u_i^n)_{\Gamma_0} \right) 
    \\ & \leq & \frac{1}{\alpha} \left( \frac{C_m}{\Delta t} ||u_i^{n-1} ||_{V_{i,0}} + ||f_i^n||_{L^2(\Omega_i)} + ||g||_{L^2(\Gamma_0)}\right) || u_i^n ||_{V_{i,0}}.
\end{eqnarray*}
This leads to the following estimate
\begin{eqnarray}
    \label{eq: u_estimate}  ||u_i^n||_{V_{i,0}} \leq \frac{1}{\alpha} \left( \frac{C_m}{\Delta t} ||u_i^{n-1} ||_{V_{i,0}} + ||f_i^n||_{L^2(\Omega_i)} + ||g||_{L^2(\Gamma_0)}\right). 
\end{eqnarray}
Similarly, by using \eqref{eq:b_infsup}, \eqref{eq:m_cont}, \eqref{eq:a_cont}, \eqref{eq:c_cont} and equations \eqref{eq:state_new}, we obtain
\begin{eqnarray*}
    ||p_i^n||_{Q_{i}} &  \leq &  \frac{1}{\beta} \sup\limits_{v_i \in V_{i,0}\backslash\{0\}} \frac{b_i(v_i, p_i^n)}{||v_i||_{V_{i,0}}} \leq \frac{1}{\beta} \sup\limits_{v_i \in V_{i,0}\backslash\{0\}} \frac{\frac{1}{\Delta t}|m_i(u_i^n - u_i^{n-1}, v_i)|   }{||v_i||_{V_{i,0}}}
  \\   & + & \frac{1}{\beta} \sup\limits_{v_i \in V_{i,0}\backslash\{0\}} \frac{|a_i(u_i^n, v_i)| + |c_i(u_i^n, u_i^n, v_i)|+|(f_i^n,v_i)_{\Omega_i}| }{||v_i||_{V_{i,0}}}
    \\ & + & \frac{1}{\beta} \sup\limits_{v_i \in V_{i,0}\backslash\{0\}} \frac{ |(g,v_i)_{\Gamma_0}| }{||v_i||_{V_{i,0}}} \leq \frac{1}{\beta}\left(\frac{C_m}{\Delta t} + C_a + C_c ||u_i^n||_{V_{i,0}}\right)||u_i^n||_{V_{i,0}}
    \\ & + & \frac{1}{\beta} \left(\frac{C_m}{\Delta t} || u_i^{n-1}||_{V_{i,0}}+||f_i^n||_{L^2(\Omega_i)} + ||g||_{L^2(\Gamma_0)} \right),
    \end{eqnarray*}
which together with the estimate \eqref{eq: u_estimate} leads to
\begin{eqnarray}
  \label{eq:p_estimation} ||p_i^n||_{Q_{i}} &  \leq & \frac{1}{\beta} \left[  \left(1 + \frac{1}{\alpha} \left(\frac{C_m}{\Delta t}+C_a\right)\right)\left(\frac{C_m}{ \Delta t} ||u_i^{n-1}||_{V_{i,0}}+||f_i^n||_{L^2(\Omega_i) } \right.  \right.
 \\ \nonumber& + & \left. \left. ||g||_{L^2(\Gamma_0) }\right) + \frac{C_c}{\alpha^2} \left( \frac{C_m}{\Delta t} ||u_i^{n-1}||_{V_{i,0}} + ||f_i^n||_{L^2(\Omega_i) } +||g||_{L^2(\Gamma_0) }\right)^2  \right] .
\end{eqnarray}
\end{itemize}

\subsection{Existence of optimal solutions}
In this subsection, we prove the existence of optimal solutions for the regularised functional \eqref{eq:functional}. The proof follows the methodology presented by Gunzburger et al. \cite{GUNZBURGER199977}.
Firstly, we define the admissibility set as follows:
\begin{eqnarray*}
    \mathcal U_{ad} = \left\{ (u_1^n ,p_1^n ,u_2^n,p_2^n, g) \in V_{1,0}\times Q_{1}\times V_{2,0}\times Q_{2}\times L^2(\Gamma_0)  \text{ such that} \quad  \quad \quad  \right.
    \\ \left. \text {equations } \eqref{eq:state_new} \text{ are satisfied and } \mathcal J_\gamma(u_1^n, u_2^n; g) < \infty \right\}.
\end{eqnarray*}

The admissibility set is clearly non--empty, since, as it was pointed out above, the restrictions to subdomains of the monolithic solution~\eqref{eq:TD_mono_equation} and its corresponding flux on the interface belong to the set. 

Let $\left\{\left(u_1^{n, (j)}, p_1^{n, (j)}, u_2^{n, (j)},p_2^{n, (j)}, g^{(j)}\right)\right\} $ be a minimizing sequence in $\mathcal U_{ad}$, i.e.,
\begin{eqnarray*}
    \lim\limits_{j \rightarrow \infty} \mathcal J_\gamma\left(u_1^{n, (j)}, u_2^{n, (j)}, g^{(j)}\right) = \inf\limits_{(u_1^n,p_1^n, u_2^n, p_2^n, g) \in \mathcal U_{ad}} \mathcal J_\gamma(u_1^n, u_2^n, g).
\end{eqnarray*}
From the definition of the admissible set $\mathcal U_{ad }$ and the functional $\mathcal J_\gamma$, it is evident that the set $\left\{ g^{(j)}\right\}$ is uniformly bounded in $L^2(\Gamma_0)$, which in turn, due to the a priori estimates~\eqref{eq: u_estimate} and \eqref{eq:p_estimation}, implies that the sets $\left\{\left(u_i^{n, (j)}\right)\right\}$ are uniformly bounded in $V_{i,0}$ and the sets $\left\{\left(p_i^{n, (j)}\right)\right\}$ are uniformly bounded in $Q_{i}$ for $i=1,2$. Thus there exists a point $\left( \hat u_1^n, \hat p_1^n, \hat u_2^n,\hat p_2^n, \hat g\right) \in \mathcal U_{ad}$ and a subsequence  $\left\{\left(u_1^{n, (j_k)},p_1^{n, (j_k)}, u_2^{n, (j_k)}, p_2^{n, (j_k)}, g^{(j_k)}\right)\right\}$  of the minimising sequence such that for $i=1,2$
\begin{eqnarray}
\label{eq: weak_u} u_i^{n, (j_k)}  \rightharpoonup  \hat u_i^n & \text{in} & V_{i,0},
\\ \label{eq: weak_p}  p_i^{n, (j_k)}  \rightharpoonup  \hat p_i^n & \text{in} & Q_{i},
\\ \label{eq: weak_g}  g^{ (j_k)}  \rightharpoonup  \hat g & \text{in} & L^2(\Gamma_0),
\\ \label{eq: strong_u}  u_i^{n, (j_k)}  \rightarrow  \hat u_i^n & \text{in} & L^2(\Omega_i),
\\ \label{eq: strong_u_int}  u_i^{n, (j_k)}|_{\Gamma_0}  \rightarrow  \hat u_i^n|_{\Gamma_0} & \text{in} & L^2(\Gamma_0).
\end{eqnarray}
The last two results are obtained by the trace theorem and compact embedding results in Sobolev spaces, see for example \cite{Evans, LionsMagenes}.

\noindent Since the forms $m_i(\cdot, \cdot)$, $a_i(\cdot,\cdot)$ and $b_i(\cdot, \cdot)$ are bilinear and continuous  by \eqref{eq: weak_u}, \eqref{eq: weak_p} and \eqref{eq: weak_g} we obtain the following convergence results:
\begin{eqnarray*}
    m_i(u_i^{n,(j_k)}, v_i) \rightarrow m_i(\hat u_i^n, v_i) & \forall v_i \in V_{i,0},
  \\  a_i(u_i^{n,(j_k)}, v_i) \rightarrow a_i(\hat u_i^n, v_i) & \forall v_i \in V_{i,0},
    \\ b_i(u_i^{n,(j_k)}, q_i) \rightarrow b_i(\hat u_i^n, q_i) & \forall q_i \in Q_{i},
    \\ b_i(v_i, p_i^{n,(j_k)}) \rightarrow b_i(v_i,\hat p_i^n) & \forall v_i \in V_{i,0},
    \\ (g^{(j_k)}, v_i)_{\Gamma_0} \rightarrow (\hat g, v_i)_{\Gamma_0}   & \forall v_i \in V_{i,0}.
\end{eqnarray*}
Concerning the trilinear form $\tilde c_i(\cdot,\cdot,\cdot)$, we exploit integration by part twice, divergence--free conditions for $u_i^{n, (j_k)}$ and $\hat u_i^n$, and the strong convergence results \eqref{eq: strong_u}--\eqref{eq: strong_u_int}. We obtain $\forall v_i \in V_{i,0}$
\begin{eqnarray*}
  \lim\limits_{k \rightarrow \infty}\frac{1}{2}\int_{\Omega_i}  (u_i^{n,(j_k)}\cdot \nabla) v_i \cdot u_i^{n,(j_k)} d\Omega  =  \frac{1}{2}\int_{\Omega_i}  (\hat u_i^{n} \cdot \nabla) v \cdot \hat 
 u_i^{n} d\Omega, 
  \end{eqnarray*}
 \begin{eqnarray*}
    \lim\limits_{k \rightarrow \infty} \frac{1}{2}\int_{\Omega_i}  (u_i^{n,(j_k)}\cdot \nabla) u_i^{n,(j_k)} \cdot v_i d\Omega  =   \lim\limits_{k \rightarrow \infty} \frac{1}{2}\int_{\Gamma_0} \left( u_i^{n,(j_k)}\cdot v_i \right) \left(u_i^{n,(j_k)} \cdot \textbf{n}_i \right) d\Gamma
  \\   -  \lim\limits_{k \rightarrow \infty} \frac{1}{2}\int_{\Omega_i}  (u_i^{n,(j_k)}\cdot \nabla) v_i \cdot u_i^{n,(j_k)} d\Omega  =   \frac{1}{2}\int_{\Gamma_0} \left( \hat u_i^{n}\cdot v_i \right) \left( \hat u_i^{n} \cdot \textbf{n}_i \right)d\Gamma
  \\  -  \frac{1}{2}\int_{\Omega_i}  (\hat u_i^{n} \cdot \nabla) v_i \cdot \hat u_i^{n,} d\Omega =  \frac{1}{2}\int_{\Omega_i}  (\hat u_i^{n} \cdot \nabla) \hat u_i^n \cdot v_i d\Omega,
  \end{eqnarray*}
  which leads to
\begin{eqnarray*}
   \lim\limits_{k \rightarrow \infty} \tilde c_i(u_i^{n,(j_k)}, u_i^{n,(j_k)}, v_i)  = \tilde c_i(\hat u_i^n, \hat u_i^n, v_i) & \forall v_i \in V_{i,0}.
   \end{eqnarray*}
These convergence results mean that the functions $\hat u_1^n, \hat p_1^n, \hat u_2^n,\hat p_2^n, \hat g$ satisfy the state equations~\eqref{eq:state_new}. We also note that the functional $\mathcal J_\gamma$ is lower--semicontinuous, i.e
\begin{eqnarray*}
     \liminf\limits_{j \rightarrow \infty} \mathcal   J_\gamma(u_1^{n,(j_k)}, u_2^{n,(j_k)}, g^{(j_k)}) \geq \mathcal   J_\gamma (\hat u_1^n,\hat u_2^n, \hat g),
\end{eqnarray*}
which implies that
\begin{eqnarray*}
     \inf\limits_{(u_1^n,p_1^n, u_2^n, p_2^n, g) \in \mathcal U_{ad}} \mathcal J_\gamma(u_1^n, u_2^n, g) = \mathcal J_\gamma (\hat u_1^n,\hat u_2^n, \hat g).
\end{eqnarray*}
Hence, we have proved the existence of optimal solutions.

\subsection{Convergence with vanishing penalty parameter}

In the previous section, we have proved the existence of optimal solutions of the regularised function $\mathcal J_\gamma$ for any $\gamma >0 $, where the parameter $\gamma$ indicates the relative importance of the two terms entering the definition of the functional. This poses an issue in our domain--decomposition setting since the optimal solution does not satisfy the coupling condition  $u_1^n|_{\Gamma_0} = u_2^n|_{\Gamma_0}$. In this section, we prove the existence of an optimal solution to the unregularised functional $\mathcal J$ with corresponds to the functional $\mathcal J_\gamma$ with $\gamma = 0$. 

Let $(u^{n,mon}, p^{n,mon})$ be a weak solution to the monolithic equations \eqref{eq:TD_mono_equation}, and for each $\gamma>0$ we denote by $(u_1^{n,\gamma}, p_1^{n,\gamma}, u_2^{n,\gamma}, p_2^{n,\gamma}, g^\gamma)$ an optimum of $\mathcal J_\gamma$  under the constraints \eqref{eq:state_new}. We define the following functions for $i=1,2$:
\begin{eqnarray*}
    u_i^{n, mon} & := & u^{n, mon}|_{\Omega_i}, \\
    p_i^{n, mon} & := & p^{n, mon}|_{\Omega_i}, \\
    g^{mon} & := & \nu \frac{\partial u_1^{n, mon}}{\partial \textbf{n}_1} - p_1^{n, mon} \textbf{n}_1 - \frac{1}{2} (u_1^{n, mon} \cdot \textbf{n}_1) u_1^{n, mon}  \quad \text{on } \Gamma_0.
\end{eqnarray*}
 Due to optimality of the point $(u_1^{n,\gamma}, p_1^{n,\gamma}, u_2^{n,\gamma}, p_2^{n,\gamma}, g^\gamma)$, we obtain that $\forall \gamma >0$
 \begin{eqnarray*}
     \mathcal J_\gamma(u_1^{n,\gamma}, p_1^{n,\gamma}, u_2^{n,\gamma}, p_2^{n,\gamma}, g^\gamma) \leq  \mathcal J_\gamma(u_1^{n,mon}, p_1^{n,mon}, u_2^{n,mon}, p_2^{n,mon}, g^{mon}), 
 \end{eqnarray*}
 which due to the definition of $u_1^{n,mon}$ and $u_2^{n,mon}$ gives us the following bound:
 \begin{eqnarray*}
     \frac{1}{2} \int_{\Gamma_0} \left| u_1^{n, \gamma} - u_2^{n, \gamma}\right|^2 d\Gamma + \frac{\gamma}{2}\int_{\Gamma_0} \left| g^\gamma \right|^2 d\Gamma \leq \frac{\gamma}{2}\int_{\Gamma_0} \left| g^{mon} \right|^2 d\Gamma \quad \forall \gamma >0. 
 \end{eqnarray*}
 The last inequality \reviewerB{tells us } that the sequence $\{ g^\gamma: \gamma >0 \}$ is bounded in $L^2(\Gamma_0)$. Following the exact same lines of arguments as in the previous section, we are able to deduce that there is a subsequence of the original sequence (we will keep the same notation for the sake of simplicity) that converges to $(u_1^{n, \ast}, p_1^{n, \ast}, u_2^{n, \ast}, p_2^{n, \ast}, g^{ \ast}) \in V_{1,0} \times Q_1\times V_{2,0} \times Q_2 \times L^2(\Gamma_0)$ in the sense of \eqref{eq: weak_u} -- \eqref{eq: strong_u_int}.
 In addition to this, the inequality above tells us that $||u_1^{n, \gamma} - u_2^{n, \gamma} ||_{L^2(\Gamma_0)} \rightarrow 0$ as $\gamma \rightarrow 0$, which in turn yields $u_1^{1, \ast} = u_2^{1, \ast}$ a.e. on $\Gamma_0$. The non--negativity of $\mathcal J$ leads to the fact that $(u_1^{n, \ast}, p_1^{n, \ast}, u_2^{n, \ast}, p_2^{n, \ast}, g^{ \ast})$ is a global minimum of $\mathcal J$. 
 Also, it is easy to see that the following functions $u^{n, \ast} \in {H_{0,\Gamma_D}^1(\Omega)}, p^{n, \ast} \in L^2(\Omega)$, defined as 
 \begin{eqnarray*}
     u^{n, \ast} :=  \begin{cases}
          u_1^{n, \ast}, \quad \text{   in } \Omega_1 \cup \Gamma_0,  \\
          u_2^{n, \ast}, \quad  \text{   in } \Omega_2 \cup \Gamma_0,
     \end{cases}
         p^{n, \ast} :=  \begin{cases}
          p_1^{n, \ast}, \quad \text{   in } \Omega_1 \cup \Gamma_0,  \\
          p_2^{n, \ast}, \quad  \text{   in } \Omega_2 \cup \Gamma_0,
     \end{cases}
 \end{eqnarray*}
satisfy the monolithic equations \eqref{eq:TD_mono_equation} in the weak sense.

\begin{remark}[Uniqueness of optimal solutions] It is well--known that the solution to the non--stationary incompressible Navier--Stokes equation in 2D is unique \cite{Richter2017}, and it can be proved that uniqueness transfers to the implicit--Euler time--discretisation scheme with a good choice of a time--step parameter (see, for instance, \cite{hairer1987solving}). This, together with the convexity of the objective functional, leads to the uniqueness of the optimal solution discussed above.
\end{remark}

\begin{remark}[Weak formulation with ``non--symmetric'' trilinear form] 
\reviewerB{The condition~\eqref{eq:trilinear_property} was crucial for analyzing the optimal-control problem. Regarding the problem formulated in the weak form~\eqref{eq:state}, numerical experiments yield similar convergence results when employing the trilinear form~\eqref{eq:trilinear_new}. We believe that this similarity arises because we have tested scenarios where the monolithic solution $(u^n, p^n)$ is well-posed, either satisfying the condition $\int_{\Gamma_N} \left| u^n\right|^2 u^n \cdot \textbf{n}\ d\Gamma \geq 0$ or considering only Dirichlet boundary conditions (with the pressure variational space comprising solely zero-mean functions to ensure uniqueness).}
\end{remark}

\subsection{DD optimal control problem gradient}\label{optimality_system}

\both{Resorting to the Lagrangian functional and the sensitivity approaches described in~\cite{prusak2022optimisationbased}, it is easy to show that the gradient of the original DD optimal control problem~\eqref{eq:functional}--~\eqref{eq:state} at a point $g \in L^2(\Gamma_0)$ is given by
\begin{equation}
  \label{eq:gradient}  \frac{d\mathcal{J}_\gamma}{dg}(u_1^n, u_2^n; g) = \gamma g + \xi_1|_{\Gamma_0} - \xi_2|_{\Gamma_0},
\end{equation}
where $\xi_1$ and $\xi_2$ are the solutions to the following adjoint problems for $i=1,2$: given $u_1^n \in V_1$ and $u_2^n \in V_2$ solutions to the state problem~\eqref{eq:state},  find $\xi \in V_{i,0}$ and $\lambda_i \in Q_i$ solving
\begin{subequations}\label{eq:adjoint}
\begin{align}
\begin{split}
\frac{1}{\Delta t }m_i(\eta_i, \xi_i) + a_i(\eta_i, \xi_i ) & +c_i \left(  \eta_i , u_i^n, \xi_i \right)  + c_i\left( u_i^n , \eta_i, \xi_i \right)  
  \\
\label{eq:adjoint1}& +  b_i ( \eta_i, \lambda_i) =  ((-1)^{i+1}\eta_i, u_1^n - u_2^n)_{\Gamma_0},
\end{split}
 & \forall \eta_i \in V_{i,0}, \\
 \label{eq:adjoint2}  & b_i (\xi_i, \mu_i)  =0  , & \forall \mu_i\in Q_i.
\end{align}
\end{subequations}}
\both{The closed--form formula for the gradient of the objective functional $\mathcal J_\gamma$ is an essential element to tackle the optimal control problem utilizing gradient--based optimisation algorithms as the one presented in~\cite{prusak2022optimisationbased} or more effective algorithms such as Broyden--Fletcher--Goldfarb--Shanno (BFGS)~\cite{bfgs_convergence, bfgs_modified} and Newton Conjugate Gradient (CG)~\cite{CG_method, CG_acceleration} algorithms. This approach leads to an iterative procedure in which both the state and the adjoint subcomponents of the DD system can be solved independently at each optimisation iteration on each subdomain and then updated through the gradient.   }
\both{In particular, in the numerical tests we will perform, we will use the limited--memory Broyden--Fletcher--Goldfarb–Shanno (L--BFGS--B) optimisation algorithm \cite{byrd1995limited}, which shows faster convergence and higher efficiency with respect to other gradient--based algorithms.}

\section{Finite Element Discretisation}
\label{high_fidelity}
In this section, we present the Finite Element spatial discretisation for the optimal control problem previously introduced. We assume to have at hand two well--defined triangulations $\mathcal T_1$ and $\mathcal T_2$ over the domains $\Omega_1$ and $\Omega_2$, respectively, and a \reviewerB{one--dimensional discretisation} $\mathcal T_0$ of the interface $\Gamma_0$.
We can then define usual Lagrangian FE spaces $V_{i,h} \subset V_{i}$, $V_{i,0, h} \subset V_{i,0}$, $Q_{i,h} \subset Q_i$, for $i=1,2,$ and $X_h \subset L^2(\Gamma_0)$ endowed with $L^2(\Gamma_0)$--norm. Since the problems at hand have a saddle--point structure, to guarantee the well--posedness of the discretised problem, we require the FE spaces to satisfy the following inf--sup conditions: there exist $c_1,  c_2, c_3,c_4 \in \mathbb R^+$ s.t. 
\begin{equation}
\label{eq:infsup}    \inf\limits_{q_{i,h} \in Q_{i,h} \backslash \{0\}} \sup\limits_{v_{i,h} \in V_{i,h} \backslash \{0\}} \frac{b_i(v_{i,h}, q_{i,h})}{||v_{i,h} ||_{V_{i,h}} ||q_{i,h} ||_{Q_{i,h}} } \geq c_i, \quad i=1,2,
\end{equation}
\begin{equation}
  \label{eq:infsup0}  \inf\limits_{q_{i,h} \in Q_{i,h} \backslash \{0\}} \sup\limits_{v_{i,h} \in V_{i,0,h} \backslash \{0\}} \frac{b_i(v_{i,h}, q_{i,h})}{||v_{i,h} ||_{V_{i,0,h}} ||q_{i,h} ||_{Q_{i,h}} } \geq c_{i+2}, \quad i=1,2 .
\end{equation}
A very common choice in this framework is to use the so--called Taylor--Hood finite element spaces, namely the Lagrange polynomial approximation of the second--order for velocity and of the first--order for pressure. We point out that the order of the polynomial space $X_h$ will not lead to big computational efforts as it is defined on the 1--dimensional curve $\Gamma_0$ \reviewerA{by the restriction of the second--order Lagrangian polynomial approximations defined by the velocity spaces. For this reason, we have chosen conformal spaces that share the degrees of freedom on $\Gamma_0$}.

Using the Galerkin projection, we can derive the following discretised optimisation problem. 
Minimise over $g_h \in X_h$ the functional
\begin{equation}
\label{eq:functional_fem}   \mathcal J_{\gamma,h} (u_{1,h}^n, u_{2,h}^n; g_h) := \frac{1}{2}\int_{\Gamma_0} \left| u_{1,h}^n - u_{2,h}^n\right|^2 d\Gamma + \frac{\gamma}{2}\int_{\Gamma_0} \left| g_h \right|^2 d\Gamma
\end{equation}
under the constraints that $u_{i, h}^n \in V_{i,h}$, $p_{i,h}^n \in Q_{i,h}$ satisfy the following variational equations for $i=1,2$
\begin{subequations}\label{eq:state_fem}
\begin{align}
\begin{split}
       \label{eq:state_fem1}   \frac{ m_i ( u_{i,h}^n - u_{i,h}^{n-1}, v_{i,h})}{\Delta t}&+ a_i(u_{i,h}^n, v_{i,h}) + c_i(u_{i,h}^n, u_{i,h}^n, v_i)     \\ 
  + b_i(v_{i,h}, p_{i,h}^n)   & =   (f_i^n, v_{i,h})_{\Omega_i}
     + \left( (-1)^{i+1} g_h, v_{i,h} \right)_{\Gamma_0}, 
\end{split}
     &\forall v_i \in V_{i,0,h},      \\
  \label{eq:state_fem2} & b_i(u_{i,h}^n, q_{i,h})   = 0, \quad \quad  &\forall q _{i,h} \in Q_{i,h},  \\ \label{eq:state_fem3}
 & \quad \quad u_i^n    =   u_{i, D, h}^n,  \quad \quad  &\text{on}  \ \Gamma_{i, D},
\end{align}
\end{subequations}
where $u_{i,D,h}^n$ is the Galerkin projection of $u_{i,D}$ onto the trace--space $V_{i,h}|_{\Gamma_{i,D}}$.
Notice that the structure of the equations \eqref{eq:state_fem} and of the functional \eqref{eq:functional_fem} is the same as the one of the continuous case. This allows us to provide the following expression of the gradient of the discretised functional \eqref{eq:functional_fem}:
\begin{equation}
 \label{eq:gradient_fem}   \frac{d\mathcal{J}_{\gamma, h}}{dg_h}(u_{1,h}^n, u_{2,h}^n; g_h) = \gamma g_h + \xi_{1,h}|_{\Gamma_0} - \xi_{2,h}|_{\Gamma_0},
\end{equation}
where $\xi_{1,h}$ and $\xi_{2,h}$ are the solutions to the discretised adjoint problem: for $i=1,2$ find $\xi_{i,h} \in V_{i,0,h}$ and $\lambda_{i,h} \in Q_{i,h}$ that satisfy 
\begin{subequations}\label{eq:adjoint_fem}
\begin{align}
\begin{split}
    \label{eq:adjoint1_fem} &\frac{m_i(\eta_{i,h}, \xi_{i,h})}{\Delta t } + a_i(\eta_{i,h}, \xi_{i,h} ) +c_i (  \eta_{i,h} , u_{i,h}^n, \xi_i )  + c_i( u_{i,h}^n , \eta_{i,h}, \xi_{i,h} ) \\ 
&+  b_i ( \eta_{i,h}, \lambda_{i,h}) =  ((-1)^{i+1}\eta_{i,h}, u_{1,h}^n - u_{2,h}^n)_{\Gamma_0},
\end{split}
\quad \forall \eta_{i,h} \in V_{i,0, h},
\\
 \label{eq:adjoint2_fem}  & \quad b_i (\xi_{i,h}, \mu_{i,h})  =0  ,  \forall \mu_{i,h}\in Q_{i,h}. 
\end{align}
\end{subequations} 
In \eqref{eq:gradient_fem}, the restriction $\cdot |_{\Gamma_0}$ is meant as an $L^2(\Gamma_0)$--projection onto space $X_h$.
We would also like to stress that at the algebraic level, the discretised minimisation problem acts only on the finite--dimensional space $\mathbb{R}^p$ of the variable $g_h$, where $p$ is the number of Finite Element degrees of freedom that belong to the interface $\Gamma_0$. 
\reviewerA{\begin{remark}[The choice of discrete space $X_h$]
    The choice for space $X_h$ might have an important impact on the convergence properties of the iterative optimisation algorithms. For the numerical tests in this work, since the objective functional and its gradient are defined by the state or adjoint solutions, our choice is to use the FE space $X_h$ defined by the degrees of freedom shared by the FE spaces $V_{1,h}$ and $V_{2,h}$.  
    Nevertheless, it is of much interest to conduct the convergence analysis of the FE DD control problem which will be the subject of future work. 
\end{remark}
}

\section{Reduced--Order Model}
\label{ROM}
\both{As highlighted in Section~\ref{introduction}, reduced--order methods are efficient tools for significant reduction of the computational time for parameter--dependent PDEs.
This section deals with the ROM for the problem obtained in the previous section, where the state equations, namely Navier--Stokes equations, are assumed to be dependent on a set of physical parameters. We rely on the classical Proper Orthogonal Decomposition (POD) technique for generating reduced spaces for each subcomponent of the problem locally. In this section, we describe two online phases based on a Galerkin projection onto the reduced spaces and on a multilayer perceptron neural network.}

\subsection{POD--Galerkin}
\label{online}
\both{We assume to have at hand local reduced spaces constructed by the classical POD compression technique~\cite{Rozza_book, prusak2022optimisationbased, benner2017model, Haasdonk} coupled with the velocity supremiser technique as described in Appendix~\ref{lifting_supremiser} to ensure the inf--sup stability of the reduced velocity--pressure spaces. For the interested reader, we describe the POD algorithm adopted in our case in Appendix~\ref{POD}}.
We can then define the reduced functions expansions 
\begin{equation*}
(u_{1,0, N}^n, p_{1, N}^n, u_{2,0, N}^n, p_{2, N}^n, g_N) \in V_{rb}^{u_1} \times V_{rb}^{p_1} \times V_{rb}^{u_2} \times V_{rb}^{p_2} \times V_{rb}^{g} 
\end{equation*}
as
\begin{align*}
u_{i, 0, N}^n   :=  \sum\limits_{k=1}^{N_{u_i}} \underline u_{i,0,k}^n\Phi_k^{u_i}, \quad 
p_{i, N}  :=   \sum\limits_{k=1}^{N_{p_i}} \underline p_{i,k}^n\Phi_k^{p_i}, \ i=1,2, \quad   
g_{ N}   :=   \sum\limits_{k=1}^{N_{g}} \underline g_{k}\Phi_k^{g} & .
\end{align*}

In the previous equations, the underlined variables denote the coefficients of the basis expansion of the reduced solution. Then, the online reduced problem reads as follows: minimise over $g_N \in V_{rb}^g$ the functional
\begin{equation}
    \mathcal J_{\gamma, N} (u_{1,N}^n, u_{2,N}^n; g_N) := \frac{1}{2}\int_{\Gamma_0} \left| u_{1,N}^n - u_{2,N}^n\right|^2 d\Gamma + \frac{\gamma}{2}\int_{\Gamma_0} \left| g_N \right|^2 d\Gamma  \label{eq:functional_rom}
\end{equation}
where $u_{1,N}^n = u_{1,0,N}^n+l_{1,N}^n$, $u_{2,N} = u_{2,0,N}^n+l_{2,N}^n$ and  $(u_{1,0,N}, p_{1,N}, u_{2,0,N}, p_{2,N})$  satisfy the following reduced equations:
\begin{subequations}\label{eq:state_rom}
\begin{eqnarray}
  \nonumber \frac{1}{\Delta t}m_i(u_{i,0,N}^n , v_{i,N}) & + &  a_i(u_{i,0,N}^n, v_{i,N})   +    c_i(u_{i,0,N}^n, u_{i,0,N}^n, v_{i,N}) 
\\ \label{eq:state_rom1}   & + &  c_i(u_{i,0,N}^n, l_{i,N}^n, v_{i,N}) + c_i(l_{i,N}^n, u_{i,0,N}^n, v_{i,N})
\\ \nonumber  & + & b_i(v_{i,N}, p_{i,N}^n)  = (f_i^n, v_{i,N})_{\Omega_i}
      +((-1)^{i+1}g_N, v_{i,N} )_{\Gamma_0}   
      \\ \nonumber & + &  \frac{1}{\Delta t}m_i(u_{i,0,N}^{n-1} , v_{i,N}) - \frac{1}{\Delta t}m_i(l_{i,N}^n , v_{i,N})
       \\\nonumber     & - &  a_i(l_{i,N}^n, v_{i,N}^n) 
      -c_i(l_{i,N},l_{i,N},v_{i,N}), \quad \forall v_{i,N} \in V_{rb}^{u_i}  ,     \\  \label{eq:state_rom2}
    b_i(u_{i,0,N}^n, q_{i,N}) &   =&   -b_i(l_{i,N}^n, q_{i,N}), \quad \forall q_{i,N} \in V_{rb}^{p_i},   
\end{eqnarray}
\end{subequations}
where $l_{i,N}^n$ is the Galerkin projection of the lifting function $l_{i,h}^n$ (see Appendix~\ref{lifting_supremiser}) for more details) to the finite dimensional vector space $V_{rb}^{u_i}$ and $i=1,2$. 

Similarly to the offline phase, we notice that the structure of the equations \eqref{eq:state_rom} and the functional \eqref{eq:functional_rom} are the same as the ones of the continuous case, so this enables us to provide the following expression of the gradient of the reduced functional \eqref{eq:functional_rom}
\begin{equation}
 \label{eq:gradient_rom}   \frac{d\mathcal{J}_{\gamma, N}}{dg_N}(u_{1,N}^n, u_{2,N}^n; g_N) = \gamma g_{N} + \xi_{1,N}|_{\Gamma_0} - \xi_{2,N}|_{\Gamma_0},
\end{equation}
where $(\xi_{1,N},\xi_{2,N})$ are the solutions to the reduced adjoint problem: find $(\xi_{1,N}, \lambda_{1,N}$, $\xi_{2,N},  \lambda_{2,N}) \in   V_{rb}^{u_1} \times V_{rb}^{p_1} \times V_{rb}^{u_2} \times V_{rb}^{p_2} $ such that it satisfies, for $i=1,2$, 
\begin{subequations}\label{eq:adjoint_rom}
\begin{align}
\label{eq:adjoint1_rom} \frac{1}{\Delta t }m_i(\eta_{i,N}, \xi_{i,N}) &+ a_i(\eta_{i,N}, \xi_{i,N} ) +c_i \left(  \eta_{i,N} , u_{i,N}^n, \xi_i \right)  + c_i\left( u_{i,N}^n , \eta_{i,N}, \xi_{i,N} \right) \\ 
\nonumber& +  b_i ( \eta_{i,N}, \lambda_{i,h}) =  ((-1)^{i+1}\eta_{i,N}, u_{1,N}^n - u_{2,N}^n)_{\Gamma_0}, \quad \forall \eta_{i,N} \in V_{i,N}^{u_i},
\\
 \label{eq:adjoint2_rom}  & \quad b_i (\xi_{i,N}, \mu_{i,N})  =0  ,  \quad \quad \quad \quad \forall \mu_{i,N}\in V_{i,N}^{p_i}. &
\end{align}
\end{subequations}
Above, the restriction $\cdot |_{\Gamma_0}$ is meant as an $L^2(\Gamma_0)$--projection onto space $V_{rb}^g$. At the algebraic level, the reduced minimisation problem can be recast in the setting of the finite--dimensional space $\mathbb{R}^{p}$, where $p$ is the number of reduced basis functions used for the control variable $g_N$ in the online phase, that is $p=N_{g}$. 
\both{We would like also to highlight the difference in the spaces we use to approximate the adjoint equations~\eqref{eq:adjoint_rom} with respect to the stationary case described in~\cite{prusak2022optimisationbased}. 
Here, instead of creating separate spaces for adjoint velocities, we used the same reduced spaces as for the state equations~\eqref{eq:state_rom} and this renders the offline stage less computationally expensive, it requires less storage and, in addition, the numerical simulations proved to be more stable.  }

\reviewerB{Here, we would like to also highlight the possibility of treating the different subdomains and control variables in the FOM or the ROM independently \cite{prusak2024optimisationbased}, obtaining various combinations of ROM/FOM models.}

\subsection{POD--NN}\label{sec:PODNN}
In this section, we would like to give a quick overview of the POD--NN method \cite{hesthaven2018non}. After the  construction of the POD reduced spaces as described in Section~\ref{POD}, the POD--NN tries to learn the map that, given the physical parameters and time, returns the reduced coefficients of the POD projection. 
To learn this map, we consider a training set of parameters  $\Xi := \lbrace (\mu_z, t_z) \rbrace _{z=1}^{KM}$ defined by a tensor product of $K$ physical parameters and $M$ time steps. 
We compute the FOM solutions for these parameters, $U_z=(u_{1,h}(\mu_z,t_z), u_{2,h}(\mu_z,t_z), p_{1,h}(\mu_z,t_z),p_{2,h}(\mu_z,t_z)) $ for $z=1,\dots, KM$ and we project the snapshots onto each reduced space to obtain the reduced coefficients, for $i=1,2$ 
$$ \underline{u}_{i,0}(\mu_z,t_z) :=\Pi_{rb,u_i}(u_{i,h}(\mu_z,t_z) ) ,\qquad \underline{p}_{i}(\mu_z,t_z) :=\Pi_{rb,p_i}(p_{i,h}(\mu_z,t_z) ) . $$
Now, for each component we consider as input the parameters $(\mu,t)$ and as output the reduced coefficients $\underline{u}_{i,0}$ or $\underline{p}_{i}$ for $i=1,2$: we therefore define $U_{u_i, \text{output}} = \{\underline{u}_{i,0}(\mu_z, t_z)\}$ and $U_{p_i, \text{output}} = \{\underline{p}_i(\mu_z, t_z)\}$. 

\noindent Once defined the input and output training sets, we build an artificial neural network (ANN) for each component $\ast \in \{ u_1, p_1,u_2,p_2\}$ that approximates $\Xi \to U_{\ast, \text{output}}$.
Then, the POD--NN reduced solutions are defined by recovering the predicted values by these ANN in the corresponding FEM space. 
Notice that this approach does not require any optimisation algorithm in the online phase, just the evaluation of the ANN at the required parameter value and time step. 

The ANN used in this algorithm is a simple dense multilayer perceptron that consists of repeated compositions of affine operations and nonlinear activation functions \cite{goodfellow2016deep}. 
The chosen architecture contains 3 hidden layers with 40, 60 and 100 neurons, respectively. 
This means that there are 4 affine mappings between the input, hidden and output layers, and at each layer, we use the hyperbolic tangent as an activation function. 
The learning of the weights and biases of the NN is optimised using the Adam algorithm \cite{kingma2014adam}, a variation of the stochastic gradient descent. 
In both test cases of the numerical result section, we used 5000 as the maximum number of optimisation iterations (epochs) and $10^{-5}$ as target for the loss functional.  

The hyperparameters are the result of a quick optimization process. 
We observed that a lower number of layers/neurons were less accurate in representing the map of interest, while more layers were too expensive to be trained in terms of necessary epochs without resulting in more accurate networks.

\section{Numerical Results}
\label{results}
We now present some numerical results obtained by applying the two--domain decomposition optimisation algorithm to the backward--facing step and the lid--driven cavity flow benchmarks. 

All the numerical simulations for the offline phase were obtained using the software FEniCS \cite{fenics}, whereas the online phase simulations were carried out using RBniCS \cite{rbnics} and EZyRB \cite{ezyrb}.

\subsection{Backward--facing step test case}
\label{backward_facing_step}
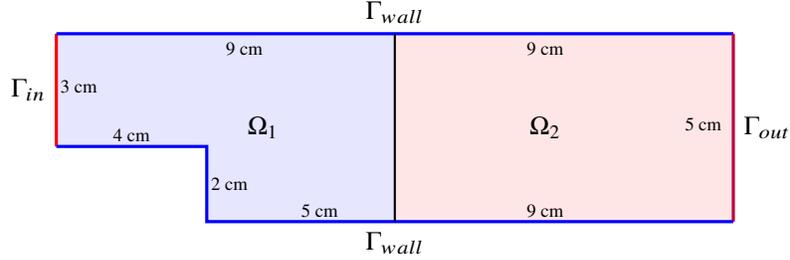
\begin{figure}
    \centering
\begin{tikzpicture}

    \fill[fill=blue!10] (0,1)  -- (0,2.5) -- (4.5,2.5) -- (4.5,0) -- (2,0)--(2,1)--(0,1); 
    \fill[fill=red!10] (4.5,0) -- (4.5,2.5)-- (9,2.5) -- (9, 1.25) -- (9,0) ;
    \draw[draw=red,very thick] (0,1) -- (0, 1.75) node [anchor=east]{$\Gamma_{in}$}  -- (0,2.5);
    \draw[draw=blue,very thick] (0,2.5) --  (4.5, 2.5) node [anchor=south] {$\Gamma_{wall}$}-- (9,2.5);
    \draw[draw=purple,very thick](9,2.5) -- (9, 1.25) node [anchor=west] {$\Gamma_{out}$}-- (9,0);
    \draw[draw=blue,very thick] (9,0) -- (4.5, 0) node [anchor=north]{$\Gamma_{wall}$} -- (2,0)--(2,1)--(0,1);
    \draw[draw=black,thick](4.5,2.5) -- (4.5, 0);

    \node[black,scale=0.7] at (1,1.15) {4 cm};
    \node[black,scale=0.7] at (2.3,0.5) {2 cm};
    \node[black,scale=0.7] at (6.5,0.15) {9 cm};
    \node[black,scale=0.7] at (3.5,0.15) {5 cm};
    \node[black,scale=0.7] at (0.3,1.8) {3 cm};
    \node[black,scale=0.7] at (2.5,2.3) {9 cm};
    \node[black,scale=0.7] at (6.5,2.3) {9 cm};
    \node[black,scale=0.7] at (8.6,1.3) {5 cm};    
    \draw (2.75, 1.25) node {$\Omega_1$};
    \draw (6.5, 1.25) node {$\Omega_2$};
\end{tikzpicture}
\caption{Physical domain and domain decomposition for the backward--facing step problem\label{fig:domain_bfs}
}
\end{figure}
We start with the backward--facing step flow test case. Figure~\ref{fig:domain_bfs} represents the physical domain of interest, the dimensional lengths and the boundary conditions. 
The splitting into two domains is performed by dissecting the domain by a vertical segment at the distance $9$ cm from the left end of the channel, as shown in Figure~\ref{fig:domain_bfs}.

We consider zero initial velocity condition, homogeneous Dirichlet boundary conditions on walls $\Gamma_{wall}$ for the fluid velocity, and homogeneous Neumann conditions on the outlet $\Gamma_{out}$, meaning that we assume free outflow on this portion of the boundary. 

We impose a parabolic profile $u_{in}$ on the inlet boundary $\Gamma_{in}$, where
$u_{in}(x,y) = 
    \left( 
         w(y),  \\
         0
    \right)^T
$
with $w(y) = \bar U \cdot \frac{4}{9} (y-2)(5-y), \ y \in [2,5]$; the range of $\bar U$ is reported in Table \ref{table:offline_bfs}.
Two physical parameters are considered: the viscosity $\nu$ and the maximal magnitude $\bar U$ of the inlet velocity profile $u_{in}$. Details of the offline stage and the finite--element discretisation are summarised in Table~\ref{table:offline_bfs}. High--fidelity solutions are obtained by carrying out the minimisation in the space of dimension equal to the number of degrees of freedom at the interface, which is 130 for our test case. 
The best performance has been achieved by using the limited--memory Broyden--Fletcher--Goldfarb–Shanno (L--BFGS--B) optimisation algorithm \cite{byrd1995limited}, where the following stopping criteria were applied: either the maximal number of iteration $\maxit$ is reached or the gradient norm of the target functional is less than the given tolerance $\tolopt$ or the relative reduction of the functional value is less than the tolerance that is automatically chosen by the \texttt{scipy} library \cite{virtanen2020scipy}. 
\begin{table}
\begin{center}
\begin{tabular}{ |r |l|| r| l| }
\hline
    \multicolumn{2}{|c||}{Physical parameters}  &  \multicolumn{2}{c|}{FE parameters} \\ \hline
    Range $\nu$ & [0.4, 2]  & Velocity--pressure space in a cell & $\mathbb P ^2 \times \mathbb P^1$\\
    Range $\bar U$ & [0.5, 4.5]  & Total dofs  & 27,890\\
    Final time $T$ & 1 & Dofs at interface & 130 \\
    & & Time step $\Delta t$ & 0.01  \\
    \hline\hline
    \multicolumn{2}{|c||}{Optimization} & \multicolumn{2}{c|}{Snapshots training set parameters}  \\ \hline
    Algorithm & L--BFGS--B & Timestep number $M$ & 100 \\
    $\maxit$ & 1000      & Parameters training set size $K$& 62\\
    $\tolopt$ & $10^{-9}$& Maximum retained modes $\nmax$ & 100 \\ \hline
\end{tabular}
\caption{Backward--facing step: computational details of the offline stage. \label{table:offline_bfs}}
\end{center}
\end{table}
\begin{figure}
    \centering
    \begin{subfigure}[b]{0.49\textwidth}
        \centering
        \includegraphics[width=0.85\textwidth]{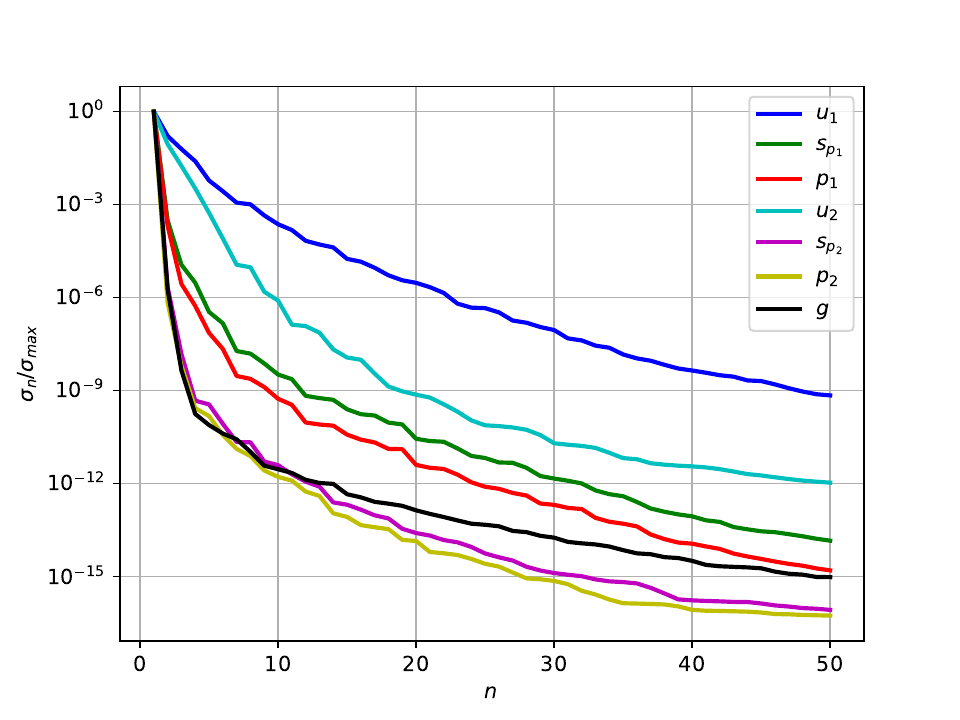}
        \caption{The singular values as a function of the number of POD modes (log scale in $y$--direction)}
         \label{fig:singlular_values_bfs}
    \end{subfigure}
    \hfill
    \begin{subfigure}[b]{0.49\textwidth}

     \begin{subfigure}[t]{\textwidth}
        \includegraphics[width=\textwidth]{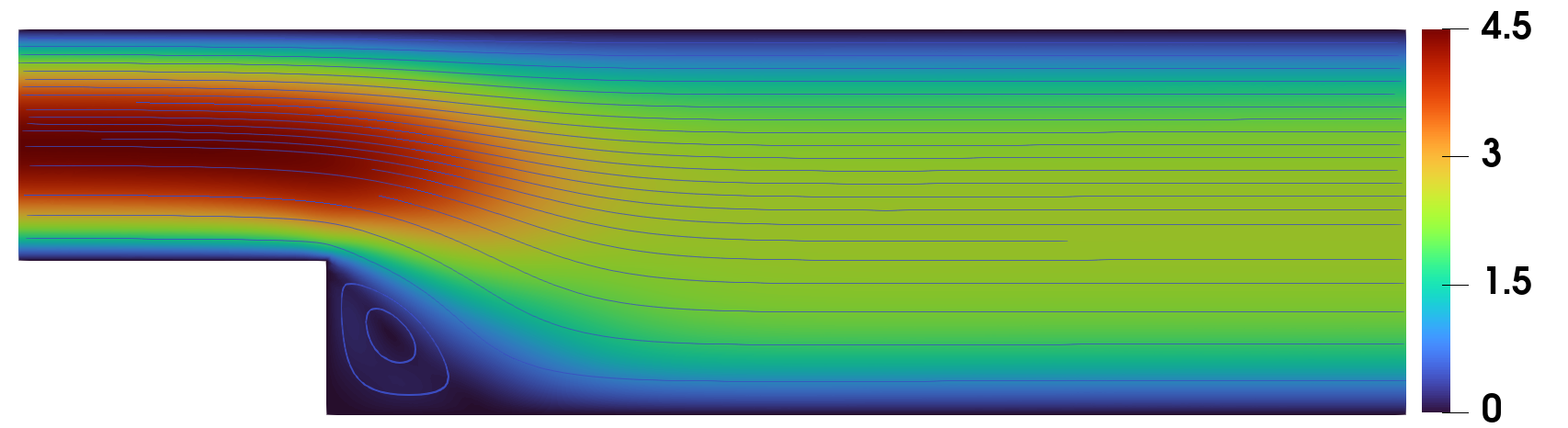}
         
    \end{subfigure}
    \vfill
     \begin{subfigure}[t]{\textwidth}
        \includegraphics[width=\textwidth]{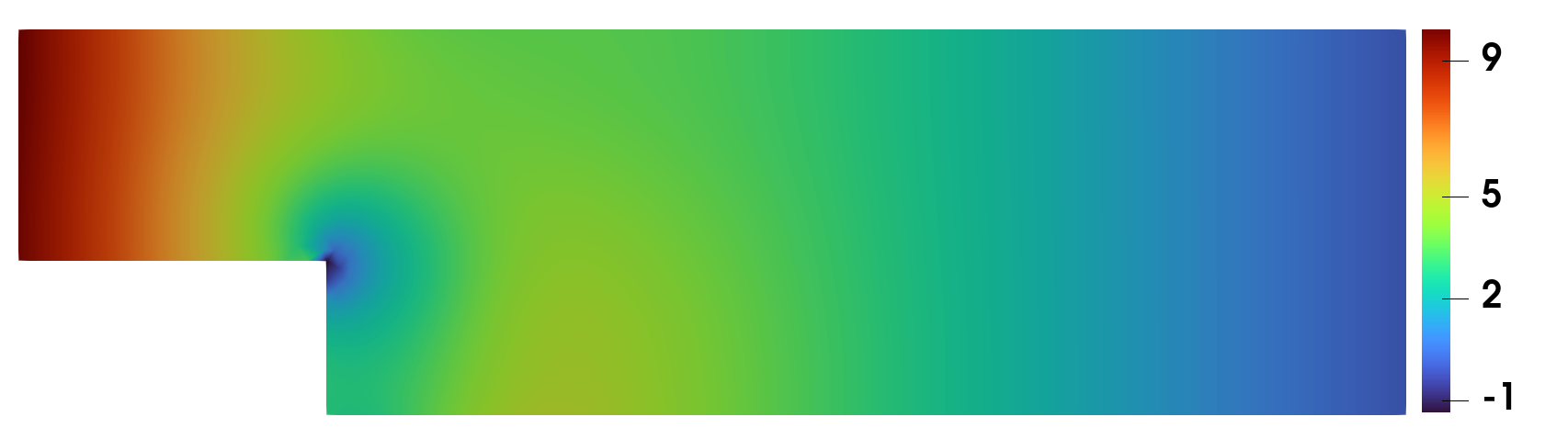}
    \end{subfigure}
    \caption{Monolithic model fluid velocity and pressure}
    \label{fig:mon_bfs}
    \end{subfigure}
    \caption{Backward--facing step: POD singular eigenvalue decay of the first 50 POD modes (a) and the monolithic solution for a parameter $(\bar U, \nu) = (4.5, 0.4)$ at the final time step (b) }
    \label{fig:pod_modes_bfs}
\end{figure}

Snapshots are sampled from a training set of $K$ parameters randomly sampled from the 2--dimensional parameter space for each time--step $t_i,\, i=1,..., M$, and the first $\nmax$ POD modes have been retained for each component. 
Figure~\ref{fig:singlular_values_bfs} shows the POD singular values for all the state and the control variables; we can see an evident exponential decay of the singular values. 
Figure~\ref{fig:mon_bfs} shows an example of a monolithic (whole--domain) solution that will be the benchmark solution for the numerical and error analysis of the DD--FOM and the ROM. 
\begin{figure}
\begin{minipage}{0.88\textwidth}
    \centering
    \begin{subfigure}[b]{0.49\textwidth}
        \includegraphics[width=\textwidth, trim={0 0 99 0}, clip]{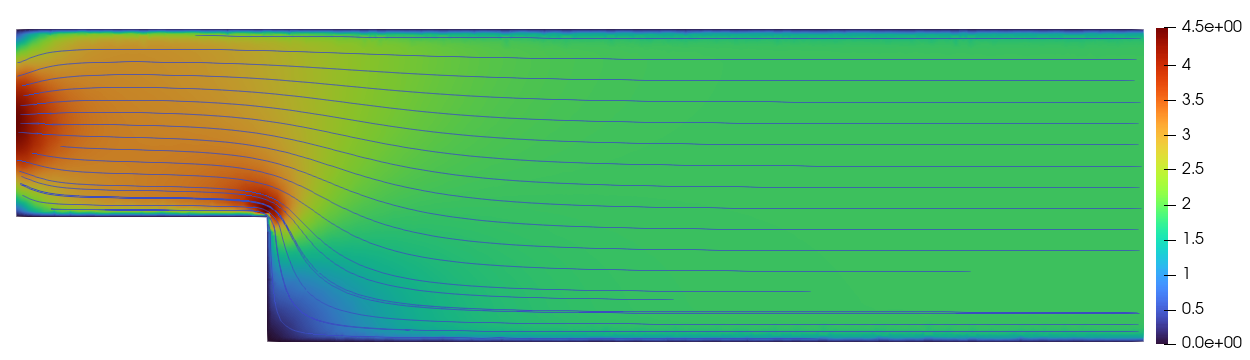}
        \caption{$t = 0.01$ }
         \label{fig:fom_bfs_u_1}
    \end{subfigure}
    \hfill
    \begin{subfigure}[b]{0.49\textwidth}
        \includegraphics[width=\textwidth, trim={0 0 99 0}, clip]{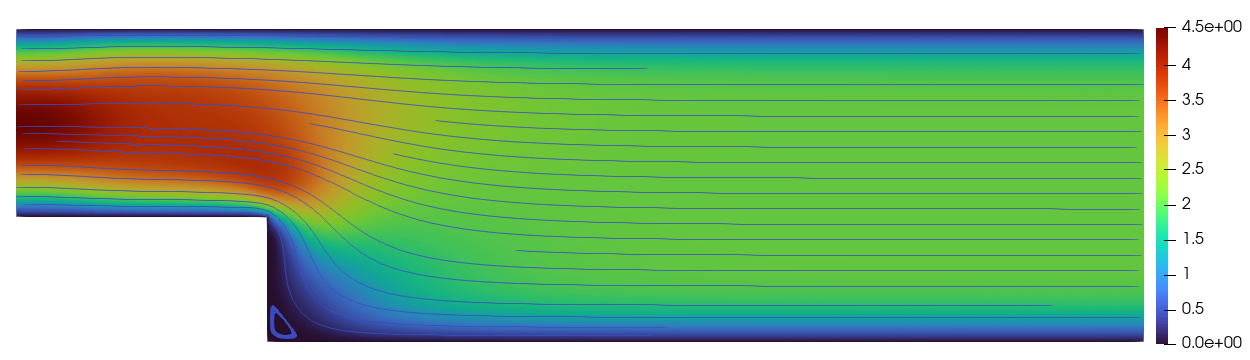}
         \caption{$t = 0.25$ }
         \label{fig:fom_bfs_u_25}
    \end{subfigure}
    \begin{subfigure}[b]{0.49\textwidth}
        \includegraphics[width=\textwidth, trim={0 0 99 0}, clip]{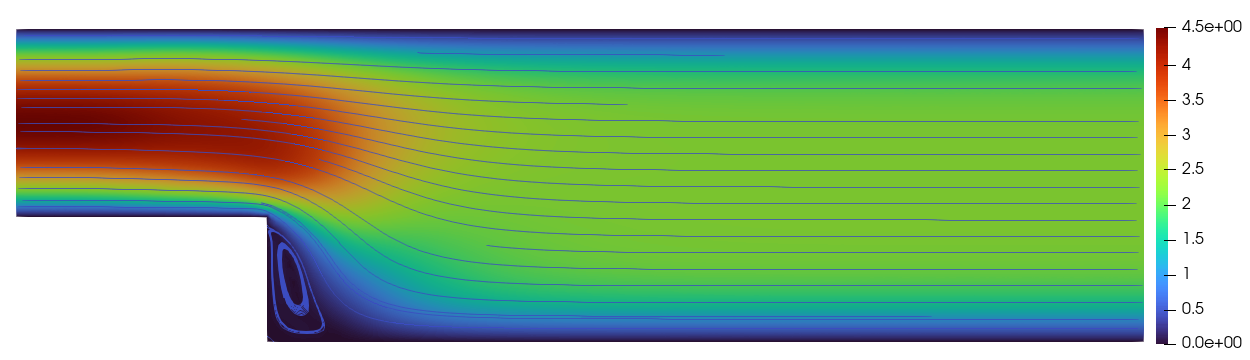}
        \caption{$t = 0.5$ }
         \label{fig:fom_bfs_u_50}
    \end{subfigure}
    \hfill
    \begin{subfigure}[b]{0.49\textwidth}
        \includegraphics[width=\textwidth, trim={0 0 99 0}, clip]{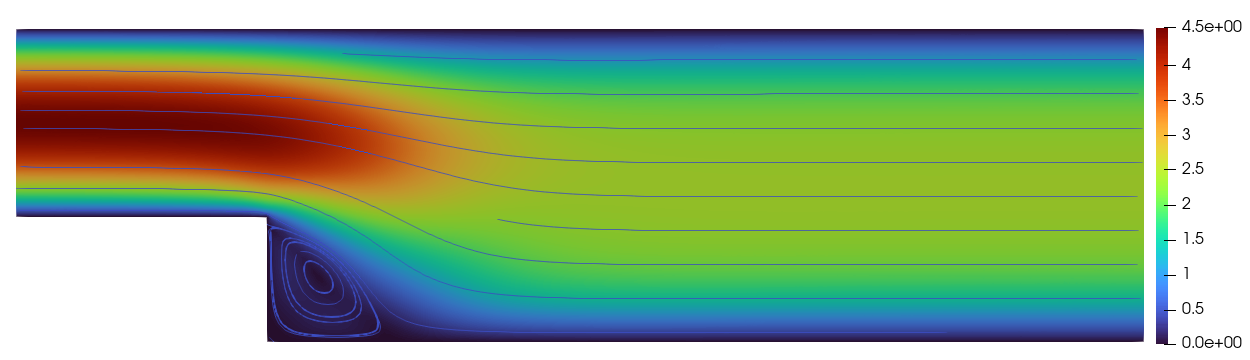}
         \caption{$t = 1$ }
         \label{fig:fom_bfs_u_100}
    \end{subfigure}
    \end{minipage}
    \begin{minipage}{0.115\textwidth}
        \includegraphics[width=\textwidth, trim={1150 0 0 0}, clip]{Figures/Truth_solution/fom_bfs_u_100.png}
    \end{minipage}
    \caption{Backward--facing step: high--fidelity solution for the velocities $u_1$ and $u_2$ at 4 different time instances}
    \label{fig:fom_bfs_u}
\end{figure}
\begin{figure}
    \centering
    \begin{subfigure}[b]{0.49\textwidth}
        \includegraphics[width=\textwidth]{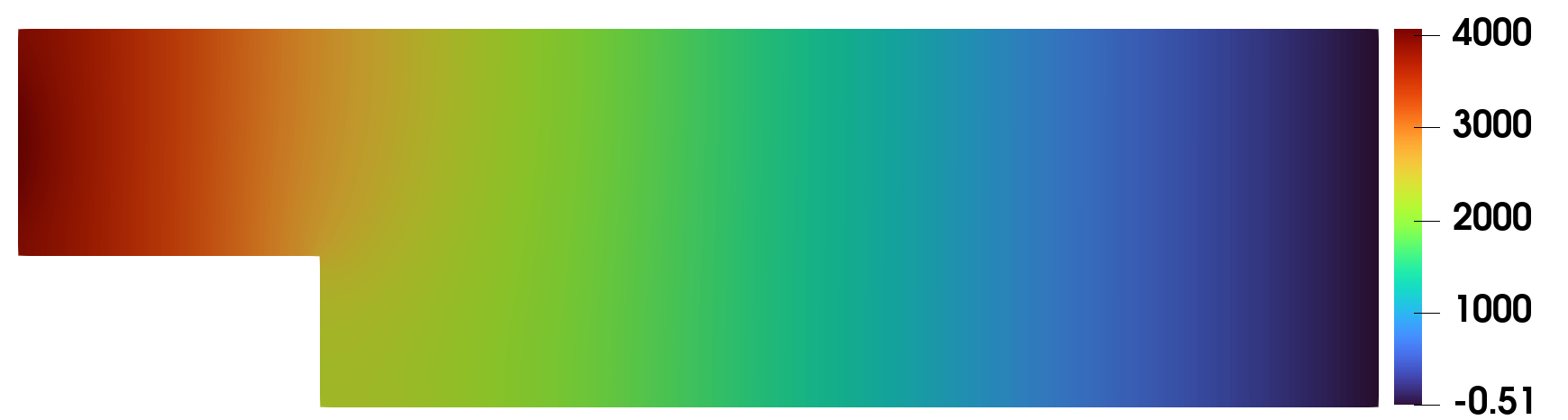}
        \caption{$t = 0.01$ }
         \label{fig:fom_bfs_p_1}
    \end{subfigure}
    \hfill
    \begin{subfigure}[b]{0.49\textwidth}
        \includegraphics[width=\textwidth]{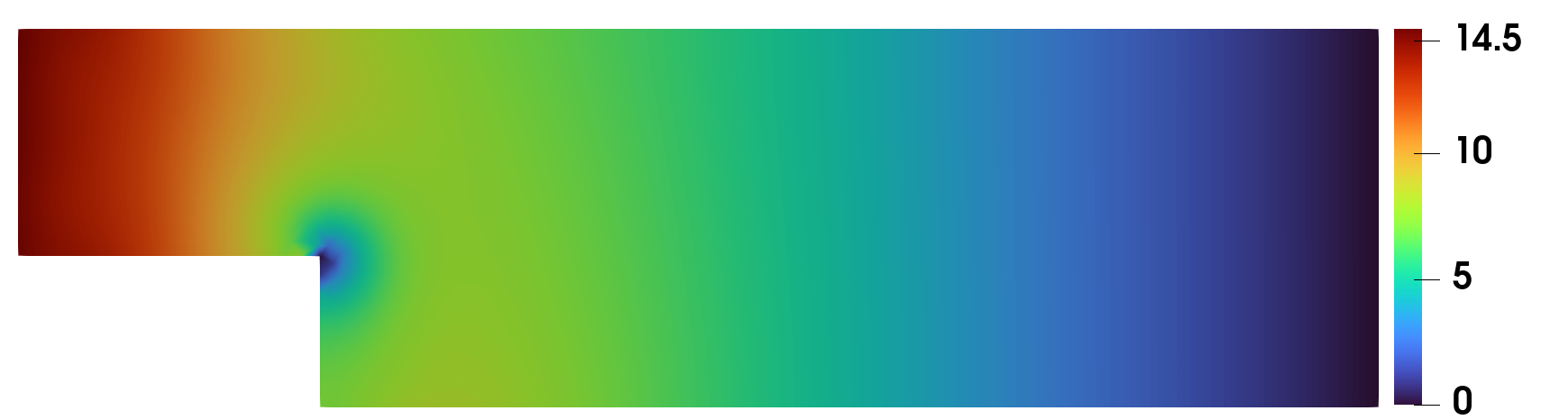}
         \caption{$t = 0.25$ }
         \label{fig:fom_bfs_p_25}        
    \end{subfigure}
    \begin{subfigure}[b]{0.49\textwidth}
        \includegraphics[width=\textwidth]{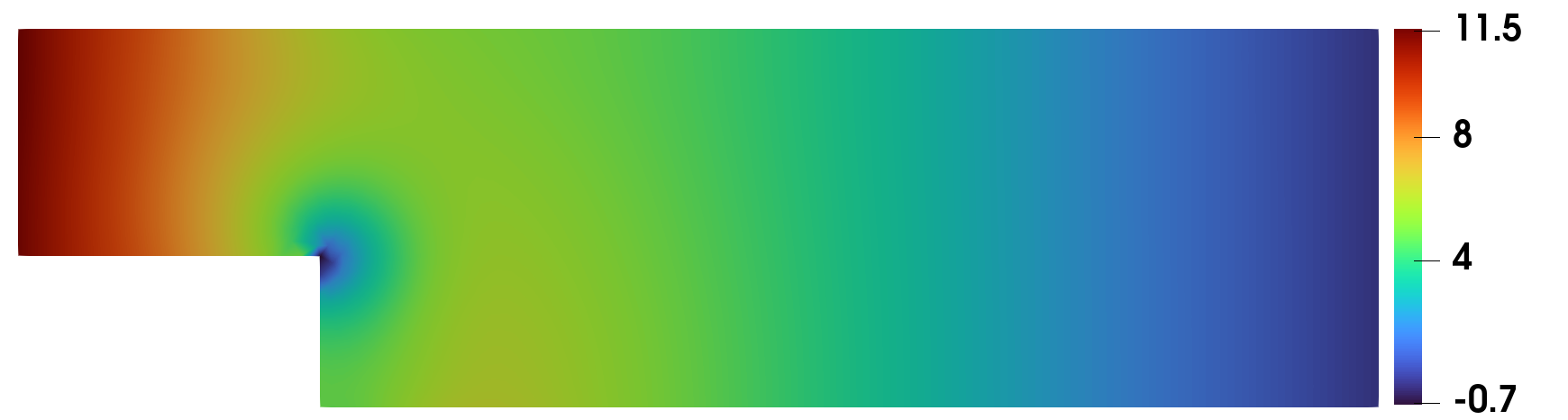}
        \caption{$t = 0.5$ }
         \label{fig:fom_bfs_p_50}
    \end{subfigure}
    \hfill
    \begin{subfigure}[b]{0.49\textwidth}
        \includegraphics[width=\textwidth]{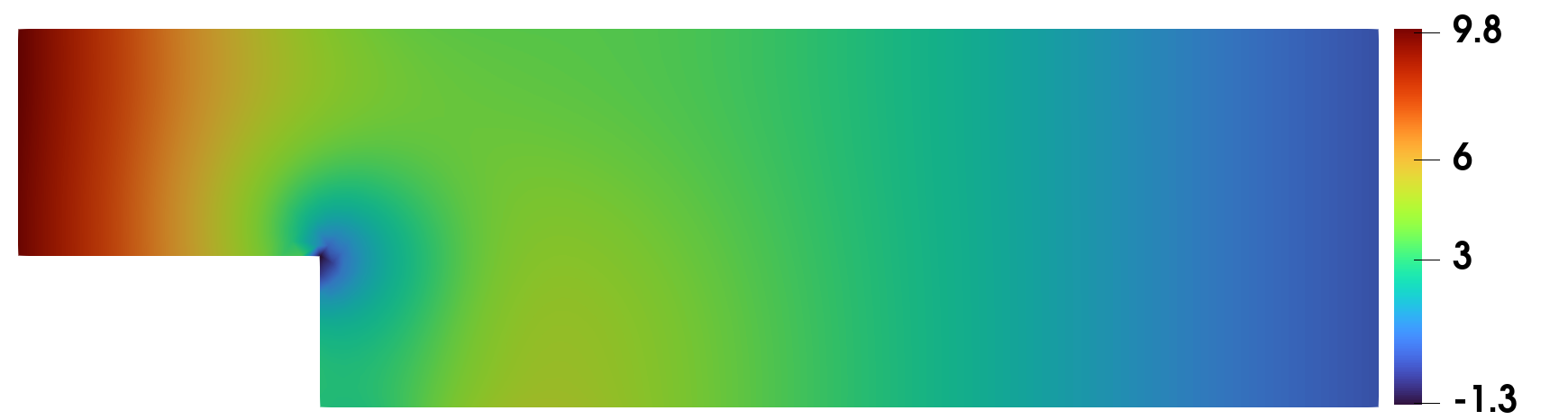}
         \caption{$t = 1$ }
         \label{fig:fom_bfs_p_100}
    \end{subfigure}
    \caption{Backward--facing step: high--fidelity solution for the pressures $p_1$ and $p_2$ at 4 different time instances}
    \label{fig:fom_bfs_p}
\end{figure}

In Table~\ref{table:online_bfs}, we list the values of the parameters for which we conduct a numerical test of the ROM and the number of POD modes for each component of the problem. The number of reduced bases is chosen so that the discarded energy for each of the components is less than $10^{-6}$. Reduced--order solutions are obtained by carrying out the minimisation in the space of dimension equal to the number of POD modes for the control $g$, which is 5 for our test case. Clearly, the minimisation in this space of dimension 5 is much simpler than in the FOM one.  The optimisation algorithm used in this test case is the same as in the FOM case described above. 

\begin{figure}
    \centering
    \begin{subfigure}[b]{0.49\textwidth}
        \includegraphics[width=\textwidth]{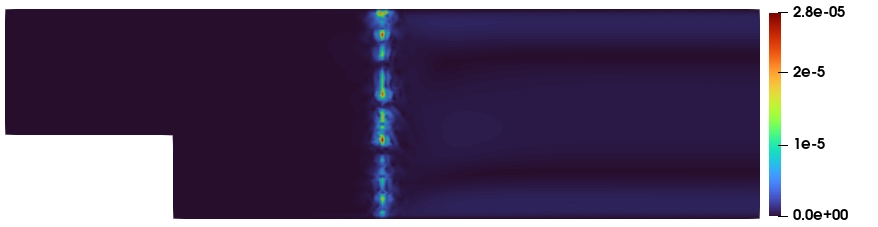}
        \caption{FOM velocity (magnitude) }
         \label{fig:fom_error_bfs_u}
    \end{subfigure}
    \hfill
    \begin{subfigure}[b]{0.49\textwidth}
        \includegraphics[width=\textwidth]{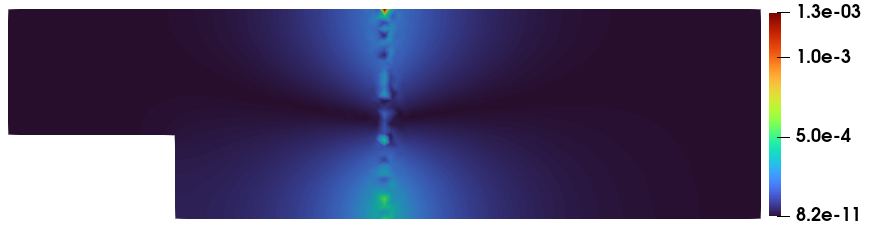}
         \caption{FOM pressure }
         \label{fig:fom_error_bfs_p}
    \end{subfigure}
    \begin{subfigure}[b]{0.49\textwidth}
        \includegraphics[width=\textwidth]{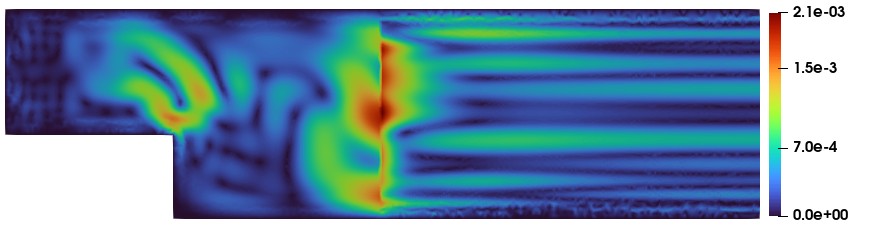}
        \caption{ROM velocity (magnitude) }
         \label{fig:rom_error_bfs_u}
    \end{subfigure}
    \hfill
    \begin{subfigure}[b]{0.49\textwidth}
        \includegraphics[width=\textwidth]{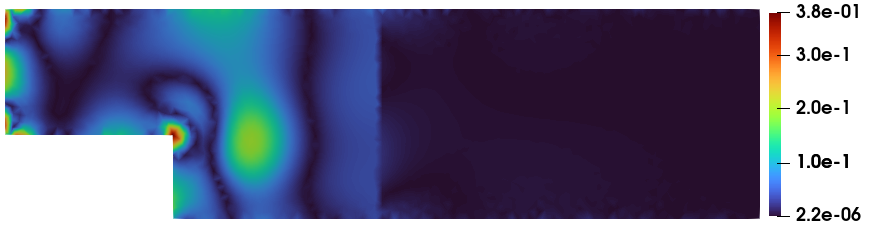}
         \caption{ROM pressure }
         \label{fig:rom_error_bfs_p}
    \end{subfigure}
    \caption{Backward--facing step: absolute errors of DD--FOM and ROM solutions w.r.t. the monolithic solution at the final time step }
    \label{fig:errors_bfs}
\end{figure}

Figures~\ref{fig:fom_bfs_u}--\ref{fig:fom_bfs_p} represent the high--fidelity solutions for a value of the parameters $(\bar U, \nu)=(4.5,0.4)$ at 4 different time instances. Visually, we can see a great degree of continuity on the interface, which will be highlighted below. 
Figure~\ref{fig:errors_bfs} shows the spatial distribution of the error with respect to the  monolithic solution at the final time step for both the FOM and ROM solutions. As expected, the error of the FOM solution is mostly concentrated at the interface, while the ROM solutions show some extra noise due to the POD reduction.
\begin{table}
\begin{center}
\begin{tabular}{ |c|c||c|c||c|c||c|c| }
\hline 
\multicolumn{2}{|c||}{Parameter}  & \multicolumn{6}{c|}{POD modes}\\ \hline \hline
 $\nu$       &  0.4   & velocity $u_1$&  30  & pressure $p_1$&  5 &  supremiser $s_1$& 5  \\ \hline
 $\bar U$  &  4.5    & velocity $u_2$& 12 & pressure $p_2$& 5  & supremiser $s_2$& 5  \\\hline
    &       & control $g $& 5  &&&& \\
    \hline
\end{tabular}
\caption{Backward--facing step: computational details of the online stage. \label{table:online_bfs}}
\end{center}
\end{table}
\begin{figure}
    \centering
        \includegraphics[width=0.7\textwidth]{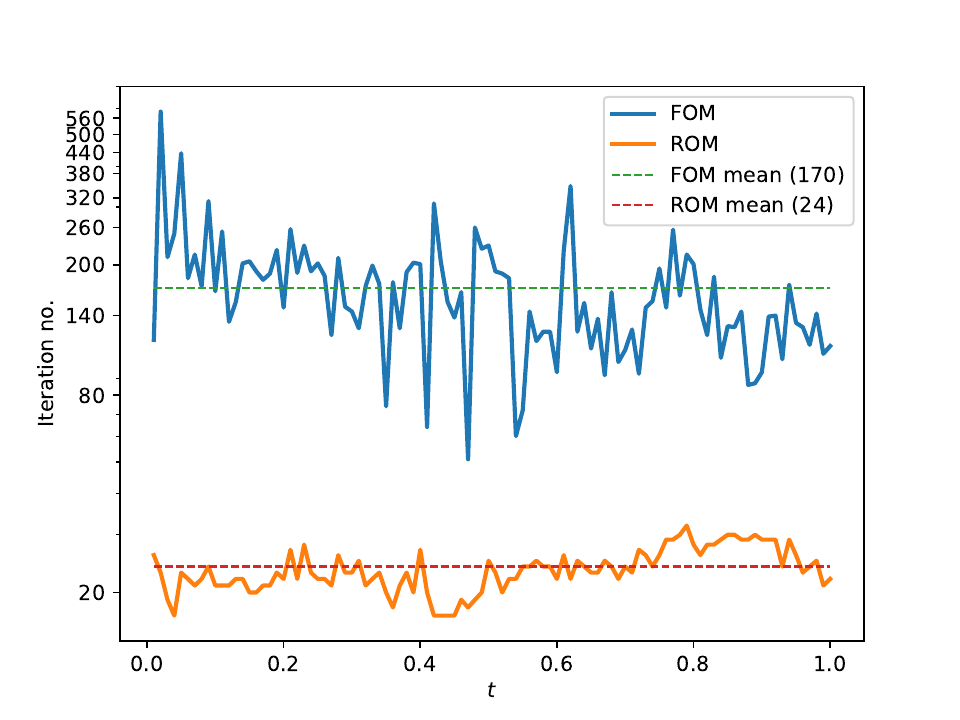}
    \caption{Backward--facing step \reviewerA{at $\nu=0.4$ and $\bar{U}=4.5$}: number of optimisation iterations of FOM and ROM solvers}
    \label{fig:its_bfs}
\end{figure}
\begin{figure}
    \centering
    \begin{subfigure}[b]{0.49\textwidth}
        \includegraphics[width=\textwidth]{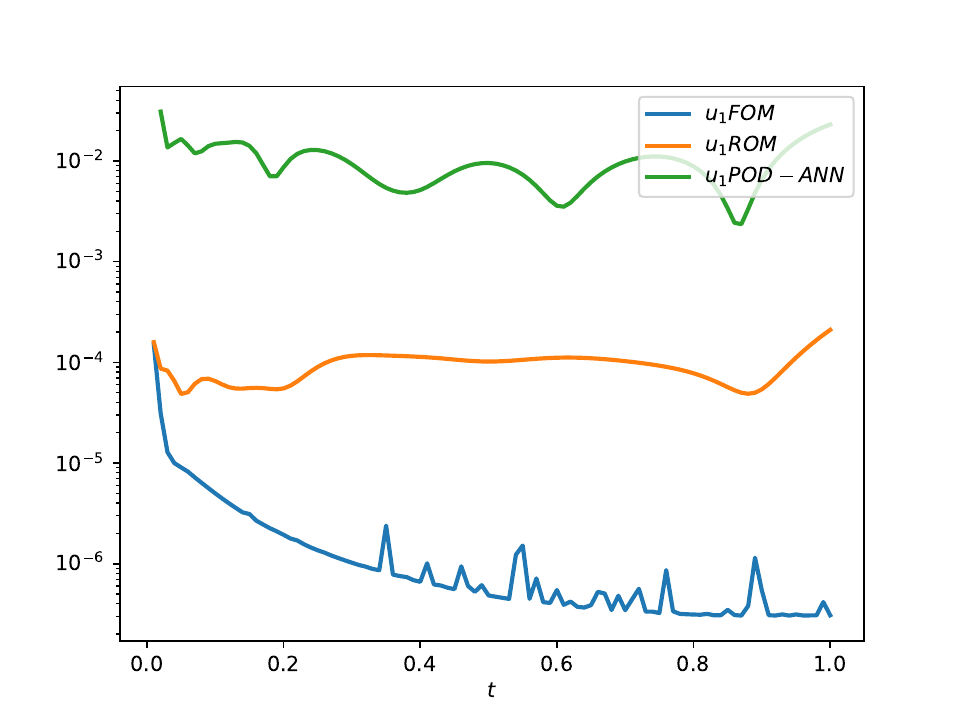}
        \caption{Velocity $u_1$ }
         \label{fig:rel_err_u1_bfs}
    \end{subfigure}
    \hfill
    \begin{subfigure}[b]{0.49\textwidth}
        \includegraphics[width=\textwidth]{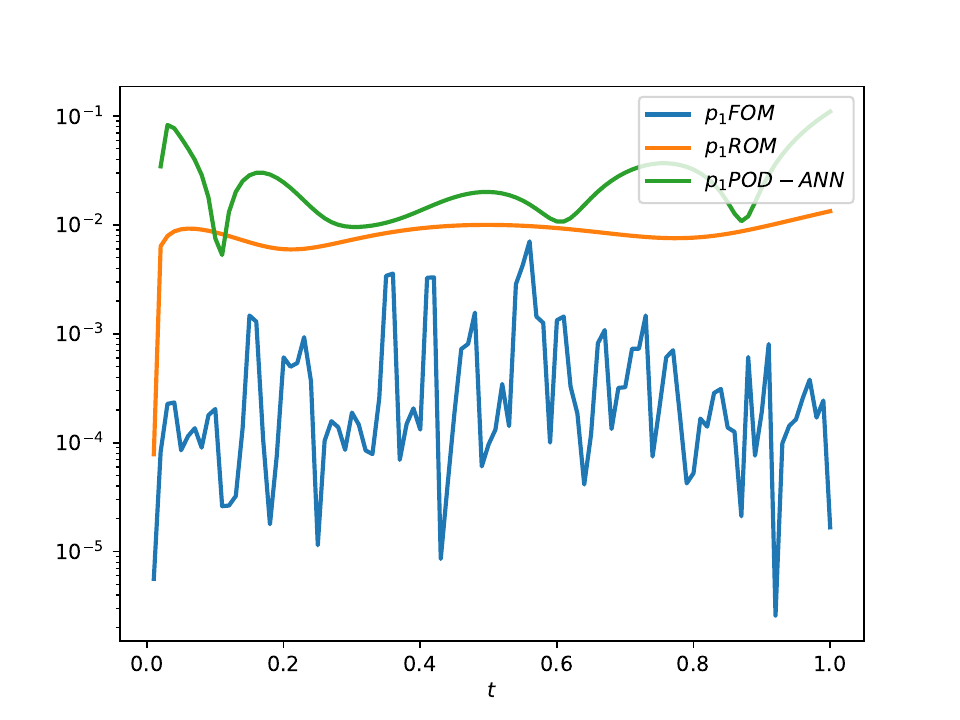}
         \caption{Pressure $p_1$ }
         \label{fig:rel_err_p1_bfs}
    \end{subfigure}
    \begin{subfigure}[b]{0.49\textwidth}
        \includegraphics[width=\textwidth]{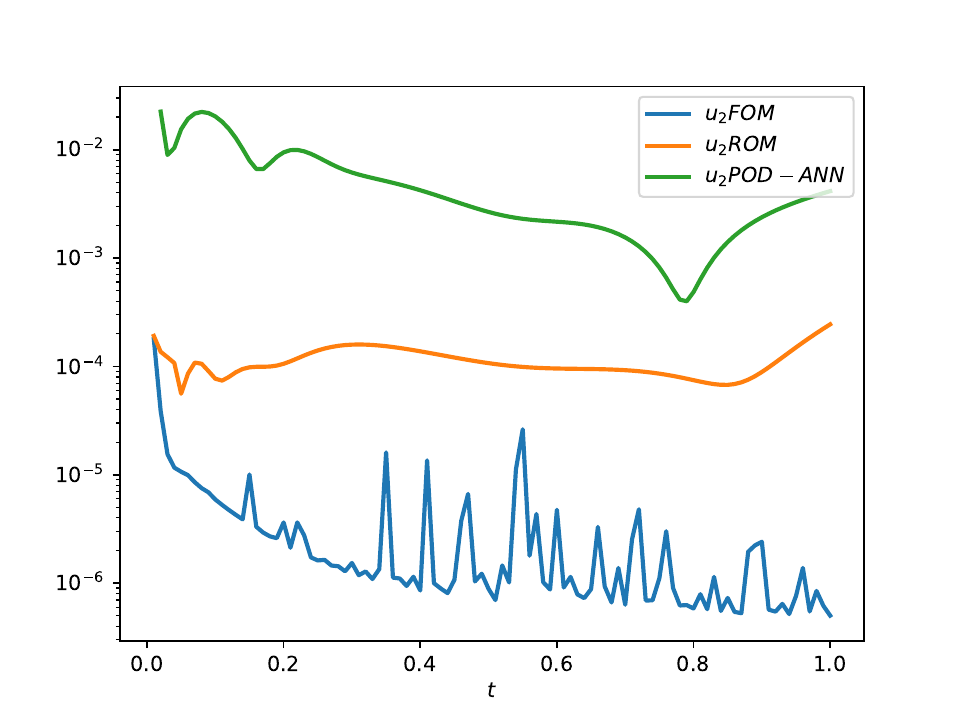}
        \caption{Velocity $u_2$ }
         \label{fig:rel_err_u2_bfs}
    \end{subfigure}
    \hfill
    \begin{subfigure}[b]{0.49\textwidth}
        \includegraphics[width=\textwidth]{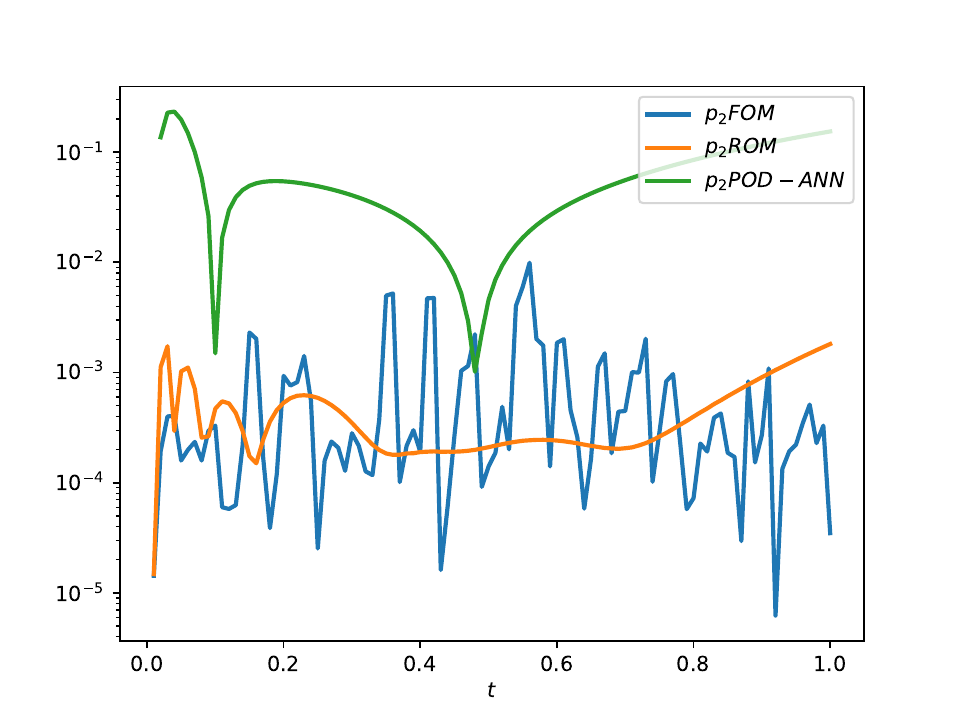}
         \caption{Pressure $p_2$ }
         \label{fig:rel_err_p2_bfs}
    \end{subfigure}
    \caption{Backward--facing step \reviewerA{at $\nu=0.4$ and $\bar{U}=4.5$}: relative errors of FOM, ROM and POD--NN solutions w.r.t. the monolithic solution  }
    \label{fig:rel_err_bfs}
\end{figure}
    
Figure~\ref{fig:its_bfs} shows the number of iterations for both FOM and ROM optimisation processes \reviewerA{for value of the parameters $(\bar U, \nu)=(4.5,0.4)$. We plot this quantity for one parameter, instead of a statistic over the parameter space, to better grasp the variability of the quantity in time, while its averages is a good representative for all other parameters}.
Each iteration of the optimisation algorithm requires at least one computation of the state and the adjoint solvers.
Therefore, we can see that we have managed to obtain a significant reduction in terms of computational efforts: the average number of iterations over all time steps in the case of the FOM solver is 170, while it is 24 in the case of the ROM solver. 
Additionally, each solver at the reduced level is of much smaller dimension (see Table \ref{table:online_bfs}), and with good use of hyper--reduction techniques (see, for instance, \cite{Rozza_book}), it will allow to obtain very efficient solvers in terms of computational time. 

We would like also to provide a comparison of the full--order and the reduced--order models with non--intrusive POD--NN model. 
Due to the discontinuity given by the initial and the boundary conditions, the first time step was excluded from the training set in order to achieve better performance. 
In practice, the first few simulations can be computed with a Galerkin projection or some FOM steps. 
Figure~\ref{fig:rel_err_bfs} shows the relative errors with respect to the monolithic solution for the FOM, ROMs and POD--NN model \reviewerA{for value of the parameters $(\bar U, \nu)=(4.5,0.4)$, again, to show the variability in time}. 
As we can see, both FOM and ROM give us very good convergence results, i.e., the relative error does not exceed 1\% in either case. 
Regarding the POD--NN, in terms of computational time, it is very effective, but the approximation can be very poor, especially in the initial and final time steps. 
Just to give an idea of the differences in the computational times, one time step of the FOM takes between 30 and 60 minutes, one time step of the ROM (without hyper--reduction) takes around 5 minutes, while a POD--NN prediction needs around $0.003$ seconds.
One of the possible scenarios could be a combination of the ROM and the POD--NN model based on the \textit{a posteriori} error estimates, so that the time steps in which the much more computationally effective ANN model fails to produce a sufficiently accurate approximation, the Galerkin projection ROM is applied. 
Similar ideas can be found \textit{inter alia} in \cite{bai2022reduced}. This will be the subject of future works. 

\subsection{Lid--driven cavity flow test case}
\label{cavity_flow}
\begin{figure}

\begin{subfigure}[t]{0.49\textwidth}
    \centering

\begin{tikzpicture}[scale=4]
      \centering
    \draw[draw=red,very thick] (0,1) -- (0.5, 1) node [anchor=south] {$\Gamma_{lid}$}  -- (1,1);
    \draw[draw=blue,very thick] (0,0) --  (0, 0.5) node [anchor=east] {$\Gamma_{wall}$}-- (0,1);
    \draw[draw=blue,very thick](0,0) -- (0.5, 0.) node [anchor=north] {$\Gamma_{wall}$}-- (1,0);
    \draw[draw=blue,very thick] (1,0) -- (1, 0.5) node [anchor=west]{$\Gamma_{wall}$} -- (1,1);

    \fill[fill=blue!10] (0,0) -- (1,0) -- (1,1) -- (0,1) -- (0,0);
    
    \draw (0.5, 0.5) node {$\Omega$};
\end{tikzpicture}
\caption{Physical domain\label{fig:Mono_domain_cavity}}
\captionsetup{justification=centering}
\end{subfigure}
\hfill        
\begin{subfigure}[t]{0.49\textwidth}
\centering    
\begin{tikzpicture}[scale=4]
  \centering
    \fill[fill=blue] (0,0) -- (1,0) -- (1,0.5) -- (0,0.5) -- (0,0);
    
    \draw (0.5, 0.25) node {$\Omega_2$};
    
    \fill[fill=red] (0,0.5) -- (1,0.5) -- (1,1) -- (0,1) -- (0,0.5);
    
    \draw (0.5, 0.75) node {$\Omega_1$};
    
    \draw (0.5, 0) [anchor=north]node {\textcolor{white}{$\Omega_2$}};

\end{tikzpicture}
\caption{Domain splitting\label{fig:dd_domain_cavity}}
\captionsetup{justification=centering}
\end{subfigure}
\caption{Lid--driven cavity flow geometry and domain decomposition}
\label{fig:geometry_cavity}
\end{figure}
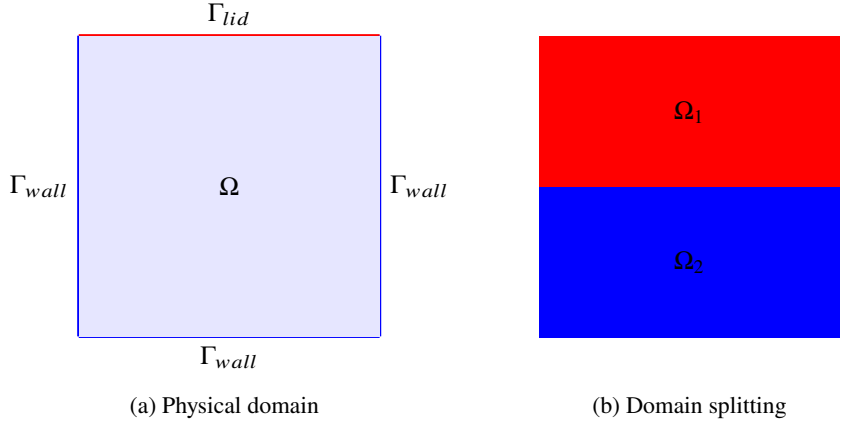
In this section, we provide the numerical simulation for the lid--driven cavity flow test case. Figure~\ref{fig:Mono_domain_cavity} represents the physical domain of interest -- the unit square. The division into two subdomains is performed by dissecting the domain by a median horizontal line, as shown in Figure~\ref{fig:dd_domain_cavity}.

We consider zero initial velocity condition, homogeneous Dirichlet boundary conditions on the boundary $\Gamma_{wall}$ for the fluid velocity and the nonzero horizontal constant velocity on the lid boundary $\Gamma_{lid}$: $u_{lid} = \left( \bar U, 0 \right)$; the values of $\bar U$ are reported in Table~\ref{table:offline_cavity}.
\begin{table}
\begin{center}
\begin{tabular}{ |r |l|| r| l| }
\hline
    \multicolumn{2}{|c||}{Physical parameters}  &  \multicolumn{2}{c|}{FE parameters} \\ \hline
    Range $\nu$ & [1, 1]  & Velocity--pressure space in a cell & $\mathbb P ^2 \times \mathbb P^1$\\
    Range $\bar U$ & [0.5, 5]  & Total dofs  & 58,056\\
    Final time $T$ & 0.4 & Dofs at interface & 294 \\
    & & Time step $\Delta t$ & 0.01  \\
    \hline\hline
    \multicolumn{2}{|c||}{Optimization} & \multicolumn{2}{c|}{Snapshots training set parameters}  \\ \hline
    Algorithm & L--BFGS--B & Timestep number $M$ & 40 \\
    $\maxit$ & 300      & Parameters training set size $K$& 10\\
    $\tolopt$ & $10^{-7}$& Maximum retained modes $\nmax$ & 100 \\ \hline
\end{tabular}
\caption{Lid--driven cavity flow: computational details of the offline stage. \label{table:offline_cavity}}
\end{center}
\end{table}
We consider one physical parameter -- the magnitude $\bar U$ of the lid velocity profile $u_{in}$. Details of the offline stage and the FE discretisation are summarised in Table~\ref{table:offline_cavity}. High--fidelity solutions are obtained by carrying out the minimisation in the space of dimension equal to the number of degrees of freedom at the interface, which is 294 for our test case. 
\begin{figure}
    \centering
    \begin{subfigure}[b]{0.49\textwidth}
        \includegraphics[width=\textwidth]{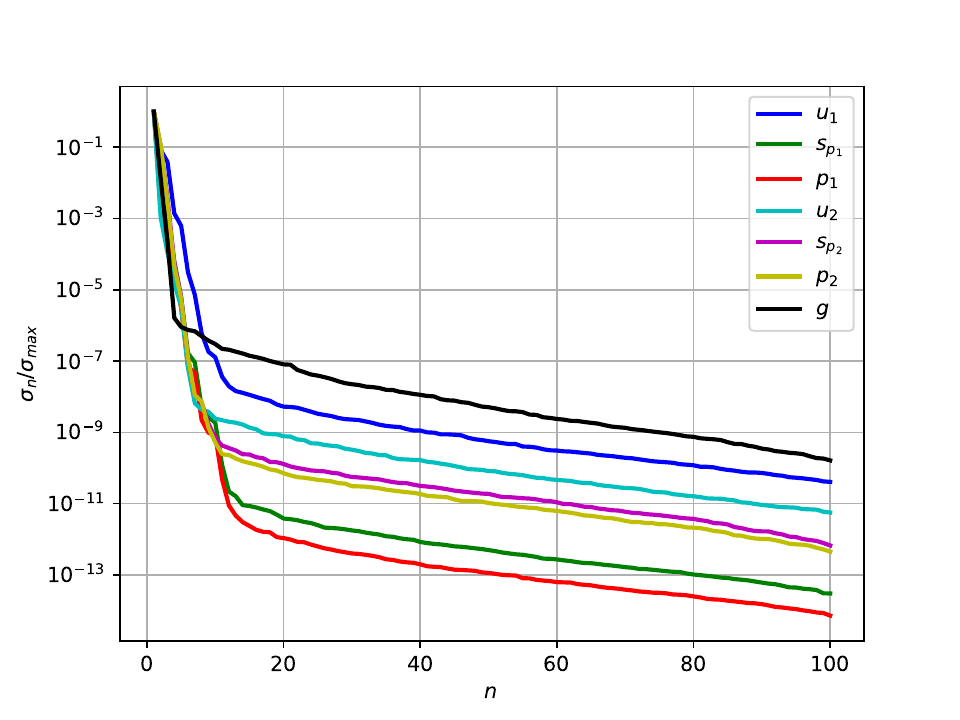}
        \caption{The singular values as a function of the number of POD modes (log scale in $y$--direction)}
         \label{fig:singlular_values_cavity}
    \end{subfigure}
    \hfill
    \begin{subfigure}[b]{0.49\textwidth}

     \quad   \includegraphics[width=0.8\textwidth]{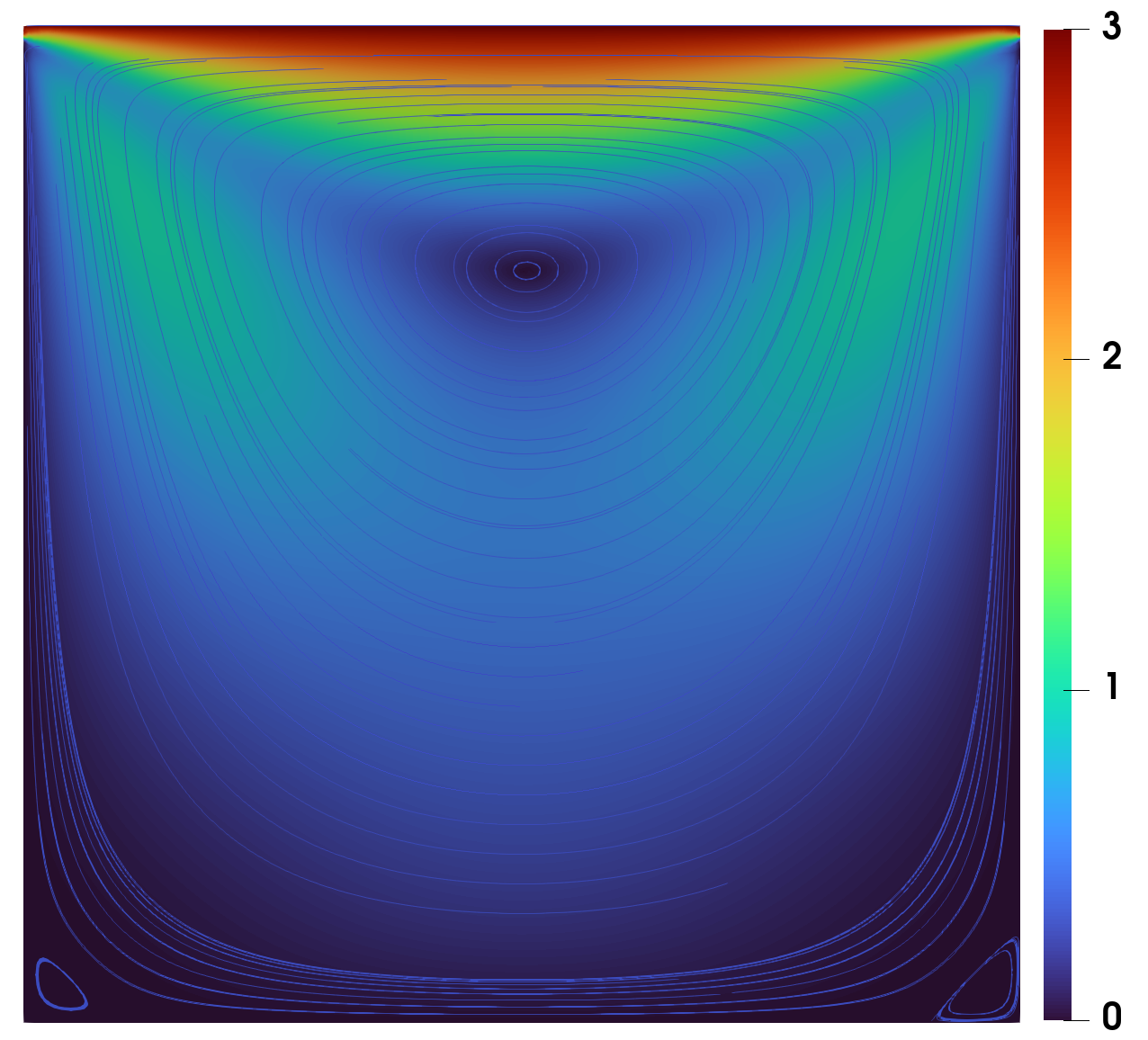}
    \caption{Monolithic model fluid velocity at the final time step  }
    \label{fig:mon_cavity}
    \end{subfigure}
    \caption{Lid--driven cavity flow: POD singular eigenvalue decay of POD modes (a) and the monolithic solution for a parameter $\bar U = 3$ at the final time step (b) }
    \label{fig:pod_modes_cavity}
\end{figure}
    \begin{figure}
    \centering
    \begin{subfigure}[b]{0.32\textwidth}
        \includegraphics[width=\textwidth]{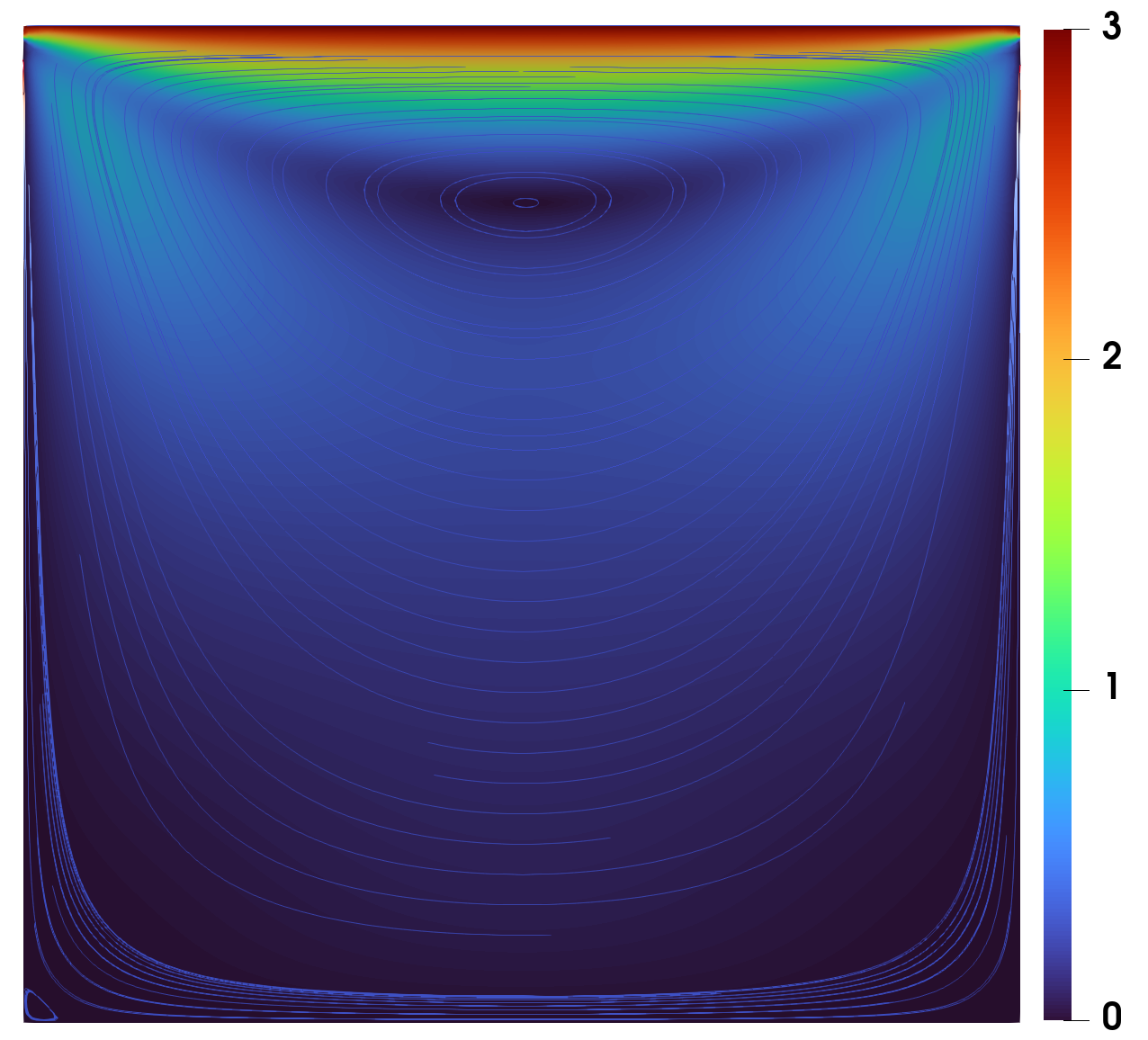}
         \caption{$t=0.01$}
         \label{fig:fom_cavity_1}
    \end{subfigure}
    \hfill
    \begin{subfigure}[b]{0.32\textwidth}
        \includegraphics[width=\textwidth]{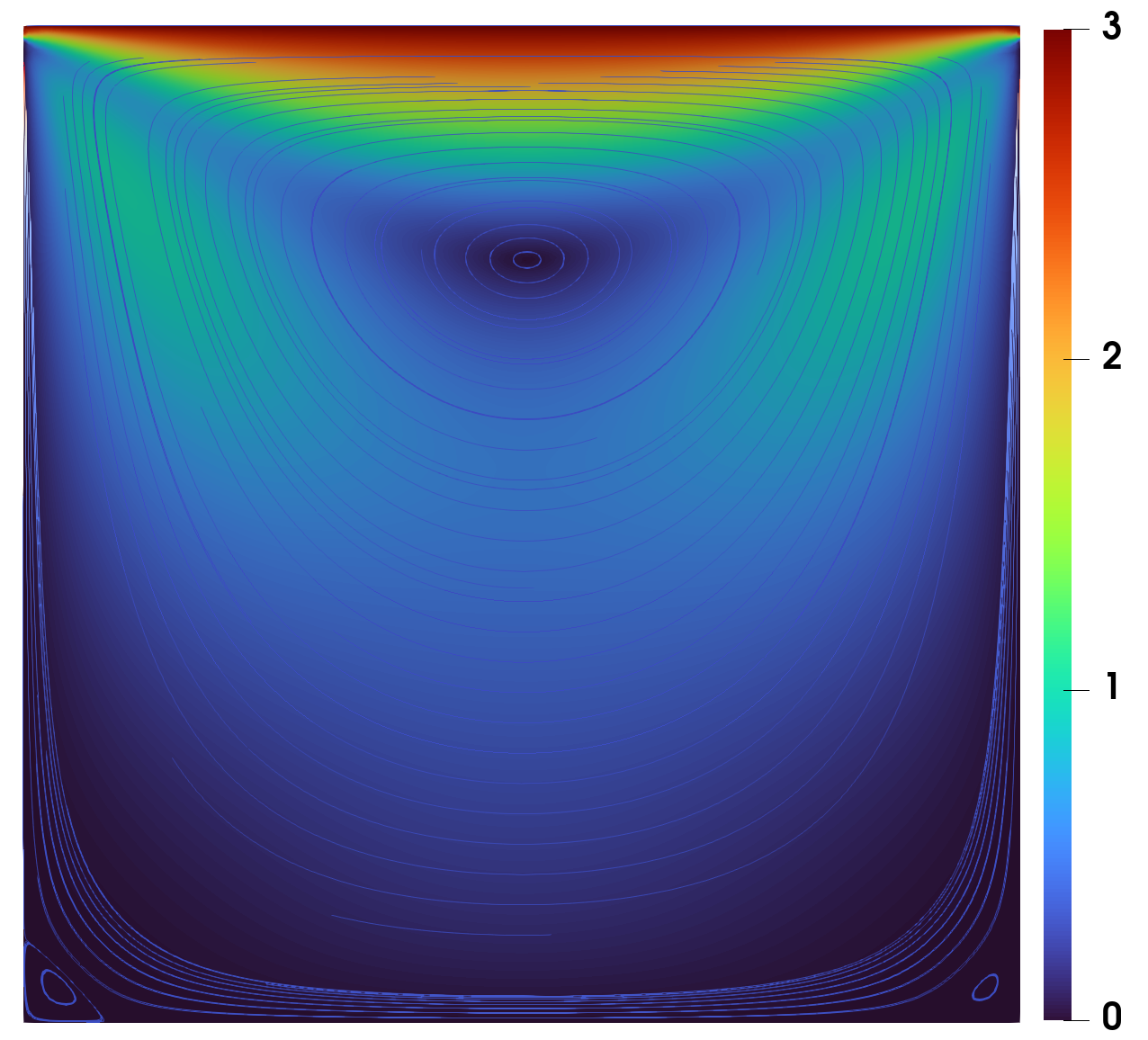}
        \caption{$t=0.05$}
         \label{fig:fom_cavity_10}
    \end{subfigure}
    \hfill
    \begin{subfigure}[b]{0.32\textwidth}
        \includegraphics[width=\textwidth]{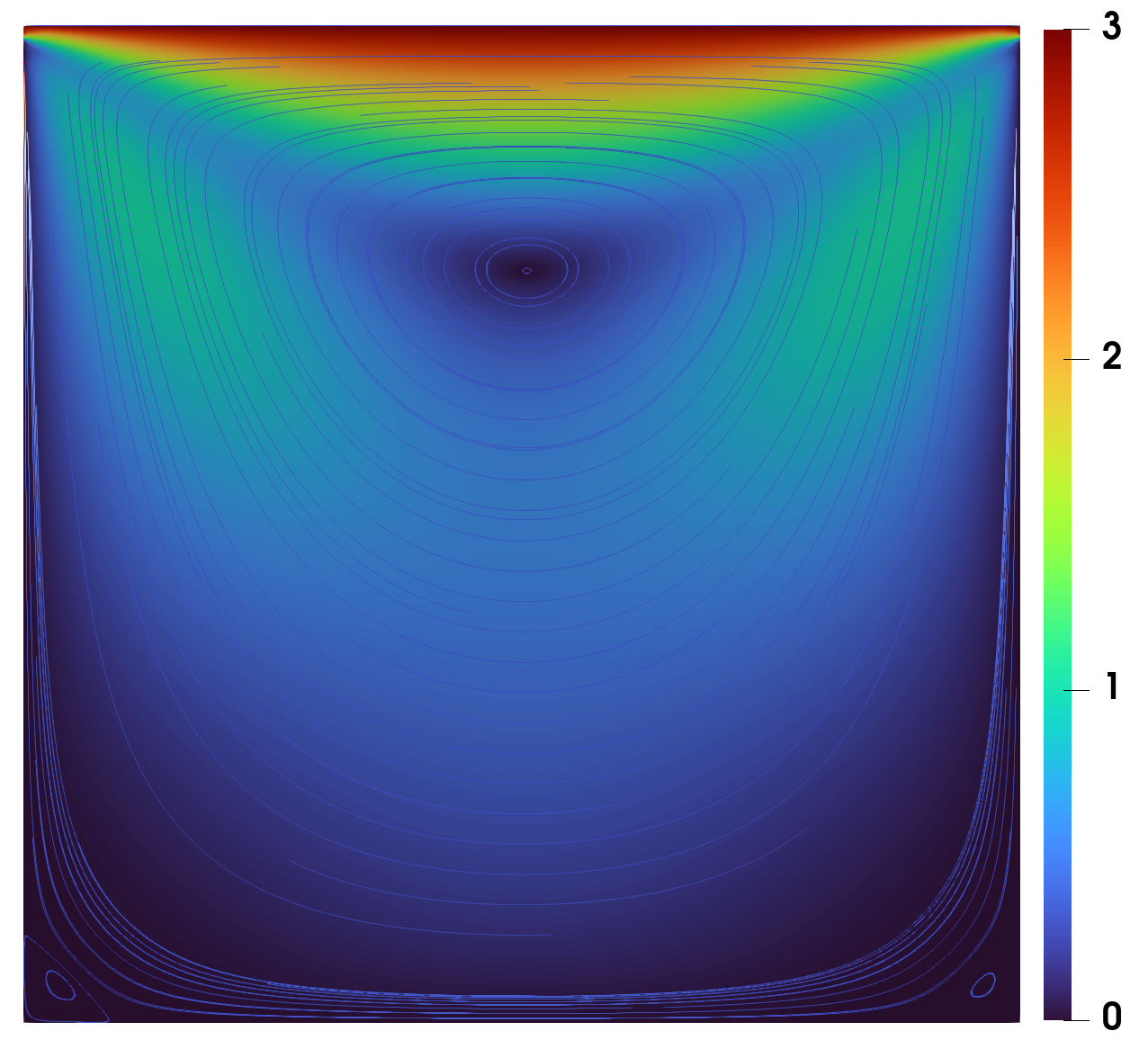}
         \caption{$t=0.4$}
         \label{fig:fom_cavity_40}
    \end{subfigure}
    \caption{Lid--driven cavity flow: FOM velocity solution at 3 different time instances}
    \label{fig:fom_cavity}
\end{figure}
Snapshots are derived from a training set of $K$ values uniformly sampled from the 1--dimensional parameter space for each time--step $t_i,\,i=1,..., M$, and the first $N_{max}$ POD modes have been retained for each component. In Figure~\ref{fig:singlular_values_cavity}, we see that the POD singular values decay even faster than in the previous test for all the state and the control variables. As before, we show in Figure~\ref{fig:mon_cavity} the monolithic (whole--domain) solution related to the parameter ($\bar U = 3$) on which we will test the DD--FOM and the ROM. 
\begin{table}
\begin{center}
\begin{tabular}{ |c|c||c|c||c|c||c|c| }
\hline 
\multicolumn{2}{|c||}{Parameter}  & \multicolumn{6}{c|}{POD modes}\\ \hline \hline
 $\nu$       &  1   & velocity $u_1$&  15  & pressure $p_1$&  10 &  supremiser $s_1$& 10  \\ \hline
 $\bar U$  &  3    & velocity $u_2$& 10 & pressure $p_2$& 10  & supremiser $s_2$& 10  \\\hline
    &       & control $g $& 5  &&&& \\
    \hline
\end{tabular}
\caption{Lid--driven cavity flow: Computational details of the online stage. \label{table:online_cavity}}
\end{center}
\end{table}
In Table~\ref{table:online_cavity}, we report the number of POD modes we use to obtain the ROM. The number of reduced bases is chosen so that the discarded energy for each of the components is less than $10^{-6}$. As before, the ROM optimization is the same used in the FOM, but on a smaller space with dimension 5 instead of 294. As optimisation algorithm, we still use the L--FBGS--B, but, in this case, we use a smaller value for $\tolopt$ of $10^{-6}$.
\begin{figure}
    \centering
        \includegraphics[width=0.85\textwidth]{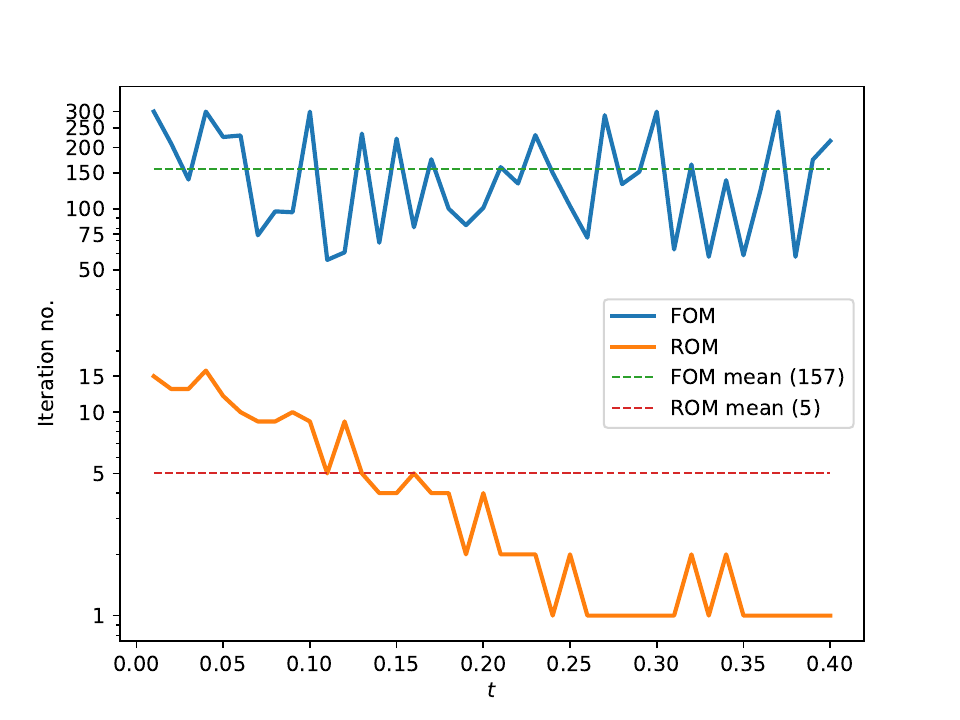}        
    \caption{Lid--driven cavity flow \reviewerA{at $\nu=1$ and $\bar{U}=3$}: number of optimisation iterations of FOM and ROM solvers }
    \label{fig:its_cavity}
\end{figure}
Figures~\ref{fig:fom_cavity} represent the DD--FOM solutions for $\bar U=3$ at 3 different time instances, where we see a qualitative agreement with the monolithic solution in Figure~\ref{fig:mon_cavity}.

Again, in Figure~\ref{fig:its_cavity} we observe that the number of optimization iterations for FOM is between 10 and 100 times larger than the ROM ones \reviewerA{for the parameter $(\bar{U},\nu)=(3,1)$}. Recalling that each iteration requires at least one computation of the state  and the adjoint solvers, we obtain a great computational advantage. For the test with $\bar U = 3$, the average number of the iteration over all time steps in the case of the FOM solver is 170 while it is 24 in the case of the ROM solver. Additionally, each solver at the reduced level is of a much smaller dimension (see Table~\ref{table:online_cavity}). 
\begin{figure}
    \centering
    \begin{subfigure}[b]{0.49\textwidth}
        \includegraphics[width=\textwidth]{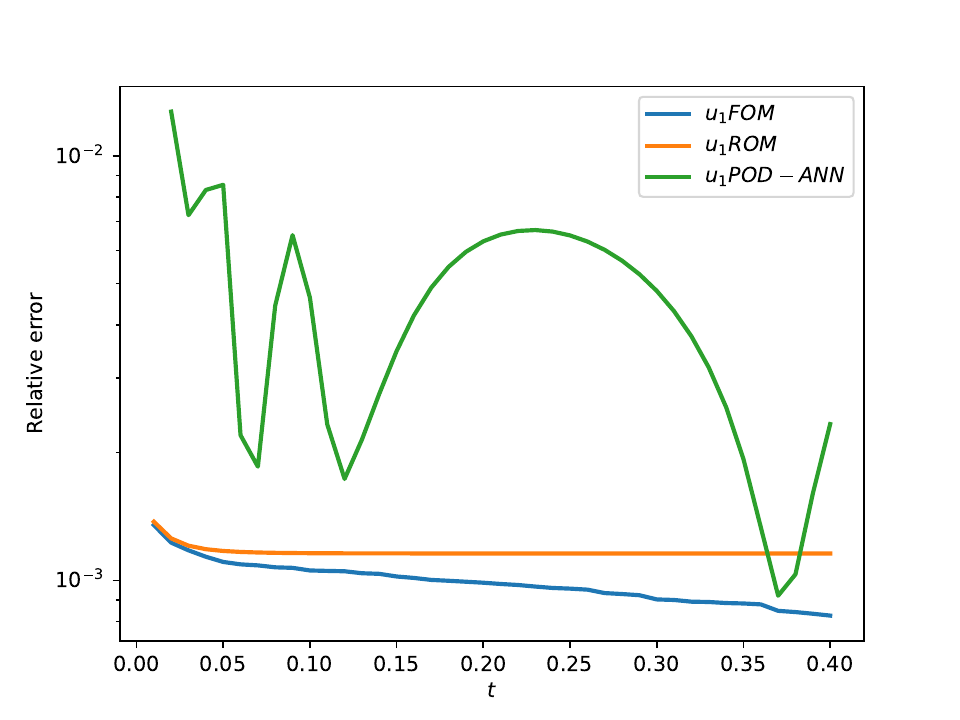}
        \caption{Velocity $u_1$ }
         \label{fig:rel_err_u1_cavity}
    \end{subfigure}
    \hfill
    \begin{subfigure}[b]{0.49\textwidth}
        \includegraphics[width=\textwidth]{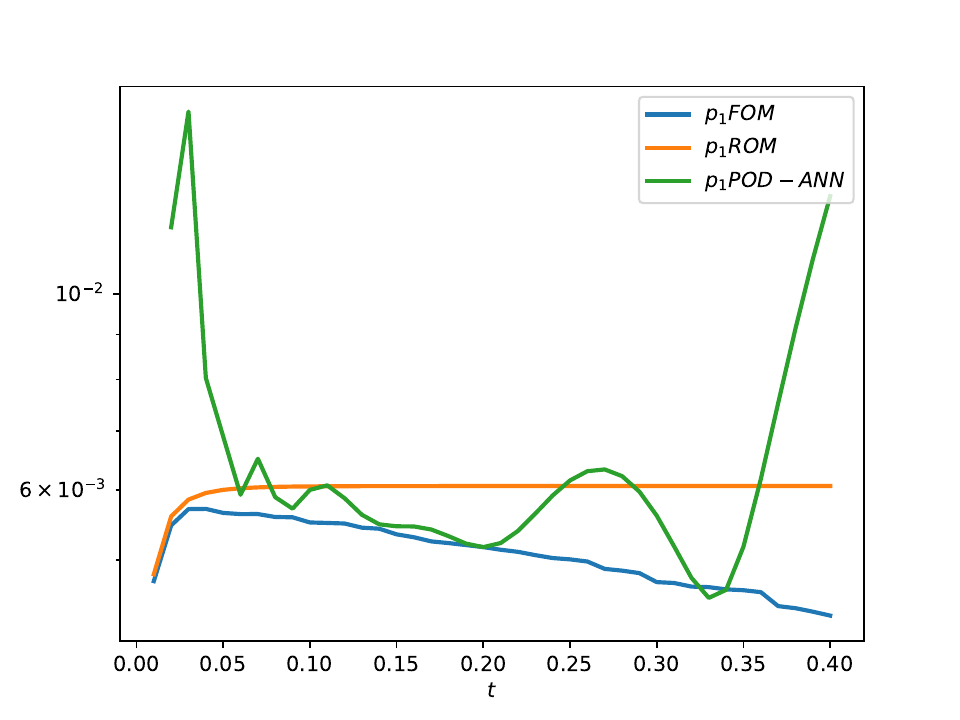}
         \caption{Pressure $p_1$ }
         \label{fig:rel_err_p1_cavity}
    \end{subfigure}
    \begin{subfigure}[b]{0.49\textwidth}
        \includegraphics[width=\textwidth]{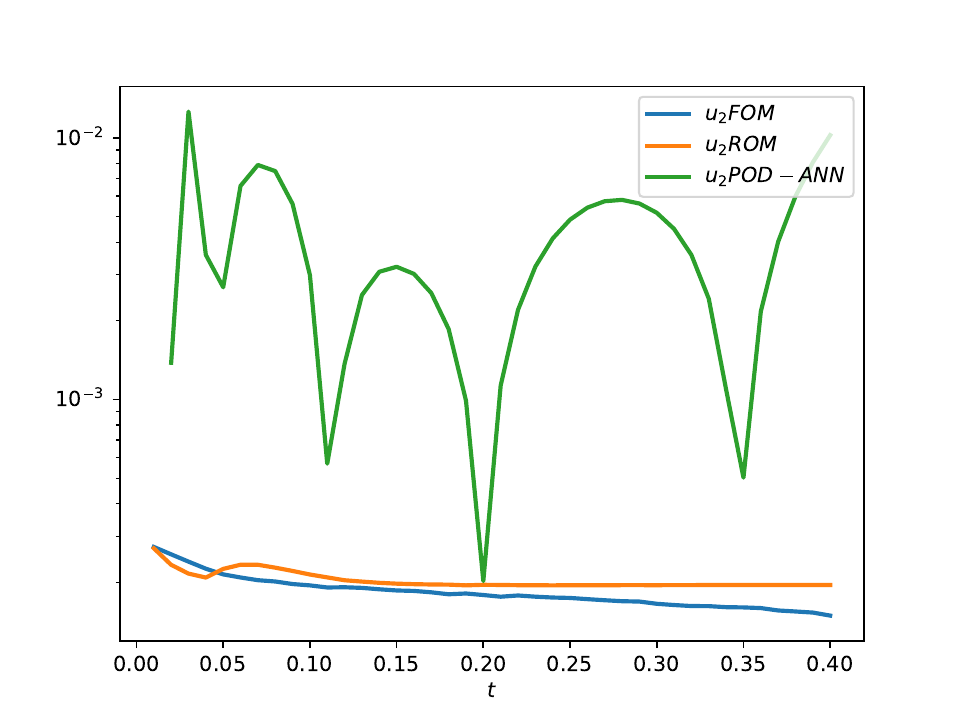}
        \caption{Velocity $u_2$ }
         \label{fig:rel_err_u2_cavity}
    \end{subfigure}
    \hfill
    \begin{subfigure}[b]{0.49\textwidth}
        \includegraphics[width=\textwidth]{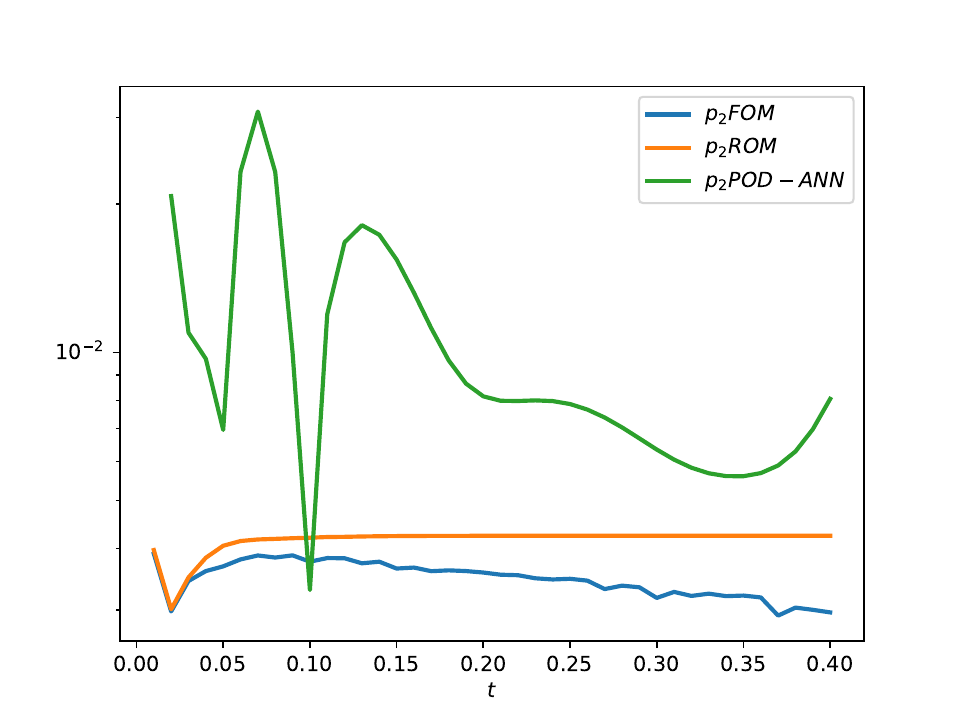}
         \caption{Pressure $p_2$ }
         \label{fig:rel_err_p2_cavity}
    \end{subfigure}       
    \caption{Lid--driven cavity flow \reviewerA{at $\nu=1$ and $\bar{U}=3$}: relative errors of FOM, ROM and POD--NN solutions w.r.t. the monolithic solution  }
    \label{fig:rel_err_cavity}
\end{figure}
As in the previous test case, we would like to provide a comparison of the full--order and the reduced--order models with non--intrusive POD--NN model. The architecture is still the one reported in Section~\ref{sec:PODNN} Again, the initial condition leads to a discontinuity in time at the starting timestep, hence, we exclude it from the training set in order to achieve better performance. Figure~\ref{fig:rel_err_cavity} shows the relative errors with respect to the monolithic solution for the FOM, ROMs and POD--NN model \reviewerA{for the parameter $(\bar{U},\nu)=(3,1)$}. As we can see, both FOM and ROM give us very good convergence results, i.e., the relative error does not exceed 1\% in either case; but, in this case, also POD--NN gives quite good results, indeed, for each variable the relative error does not exceed 3\%.  Computational times for each method, FOM, ROM and POD--NN, are comparable with those of the backward--facing step flow, that is, one time step of the FOM takes between 15 and 45 minutes, one time step of the ROM (without hyper--reduction) takes on average 5 minutes, while a POD--NN prediction needs around $0.003$ seconds.

\reviewerB{
\begin{remark}[Sensitivity to the domain splitting]
    Similarly to the stationary case analyzed in~\cite{prusak2022optimisationbased}, we encounter the issue of choosing the optimal domain splitting.  Also in the non--stationary lid--driven cavity flow case, we observe a significantly slower convergence when the domain is split by a vertical segment instead of a horizontal one. 
\end{remark}
}

\reviewerA{
\subsection{Validation of the POD--NN method}\label{sec:validation_podnn}
As we have seen in the numerical experiments above, unlike the POD--Galerkin, the POD--NN ROM seems to be very sensitive to the presence of discontinuities in time. 
In this section, we would like to investigate under which circumstances the POD--NN technique is able to produce a good approximation by considering various scenarios and performing a statistical analysis in terms of physical parameters. 
The first necessary step in all the approaches below is to provide a larger set of snapshots for generating RB spaces and for the training artificial neural network, as is described in Section~\ref{sec:PODNN}.
\begin{figure}
    \centering
    \begin{subfigure}[b]{0.49\textwidth}
        \includegraphics[width=\textwidth, trim = {26 0 24 0}, clip]{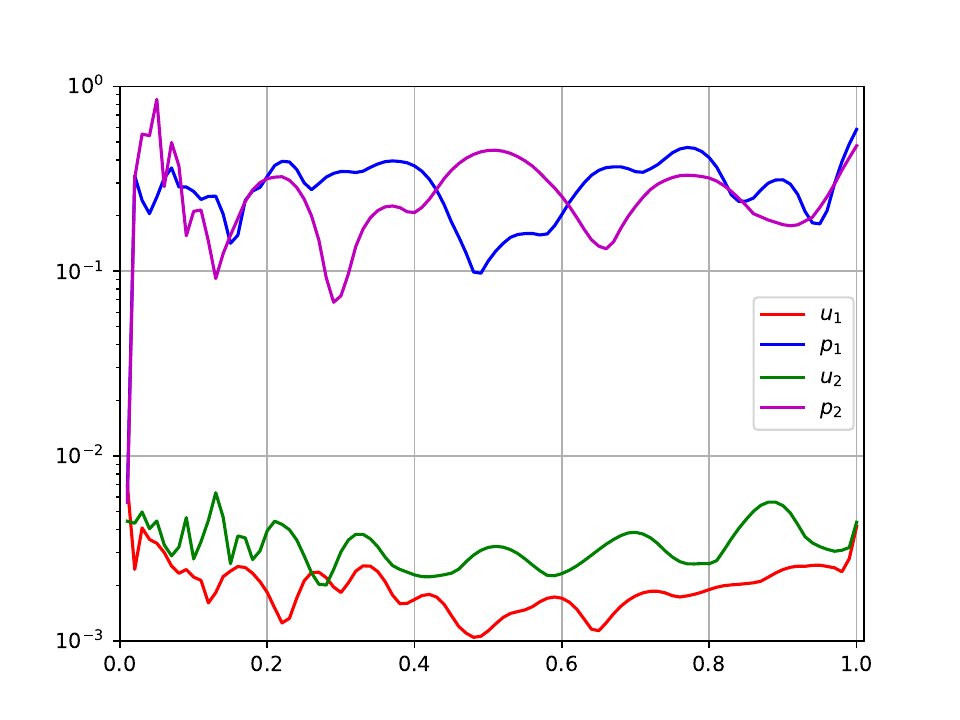}
        \caption{Mean relative errors over the entire time interval $[0,1]$}
         \label{fig:bfs_vacuum_mean_full}
    \end{subfigure}
    \hfill
    \begin{subfigure}[b]{0.49\textwidth}
        \includegraphics[width=\textwidth, trim = {15 0 35 0}, clip]{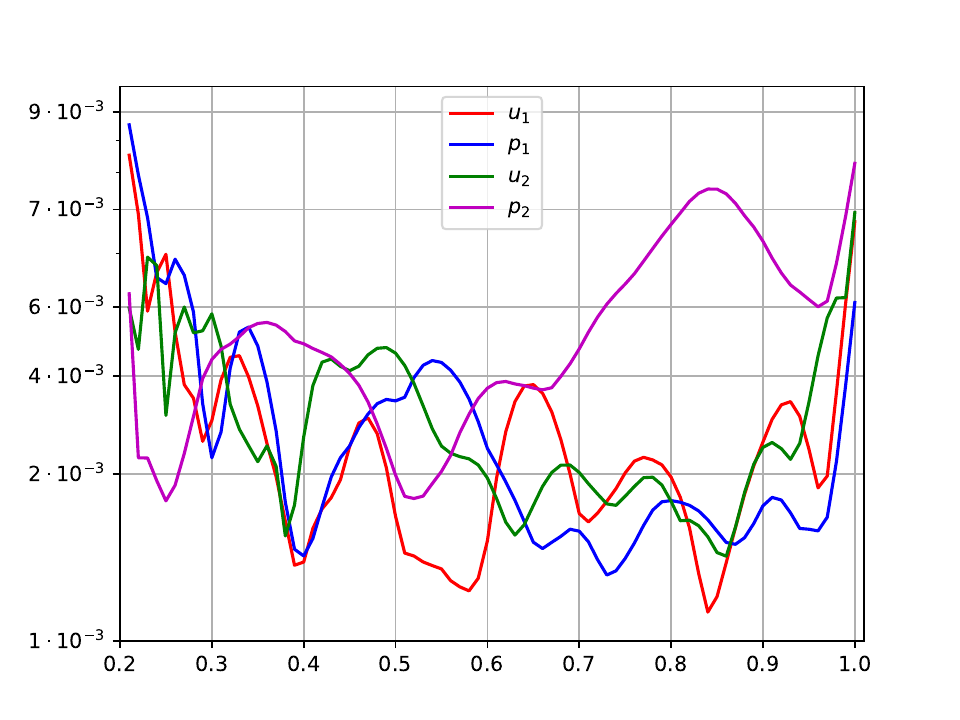}
         \caption{Mean relative errors over the restricted time interval $[0.21,1]$}
         \label{fig:bfs_vacuum_mean_restricted}
    \end{subfigure} 
    \caption{Mean over the parameter testing set of  POD--NN relative errors for the discontinuous in time backward--facing step test case  with respect to $\bar U$ over the entire (a) and restricted (b) time intervals}
    \label{fig:bfs_vacuum_mean}
\end{figure}

Regarding the backward--facing step test case presented in Section~\ref{backward_facing_step}, we have seen that the constructed NN--ROM provides very poor approximations, most probably due to the presence of a discontinuity between the inlet boundary and the initial conditions. We restrict the analysis by considering only one physical parameter, $\bar U$ to generate new POD--based reduced spaces and to train the ANN. 
In this case, we draw $101$ parameters uniformly sampled from the parameter space reported in Table~\ref{table:offline_bfs}, we randomly shuffle the parameter set and divide it into a training set $\mathbb P_{\text{train}}$, containing $75$ parameters, and a testing set $\mathbb P_{\text{test}}$, composed of remaining $26$ parameters.

We consider two scenarios: in the first one we construct the parameter space for POD compression and the ANN training by using the training parameter set and over all the time steps, i.e, the set $\mathcal P_{\text{train}}^{\text{full}} = \left\{ (\bar U, t^n): \bar U \in \mathbb P_{\text{train}},  \ n \in \{1,...,100\}\right\}$ whereas in the other one we discard the first $20$ time steps and consider only the training parameter set and the time steps $t^n$ for $n\geq 21$, i.e., the set $\mathcal P_{\text{train}}^{\text{restricted}} = \left\{ (\bar U, t^n): \bar U \in \mathbb P_{\text{train}}, \ n \in \{21,...,100\}\right\}$.
We perform the POD compression using the snapshots corresponding to parameter values in $\mathcal P_{\text{train}}^{\ast}$ for $\ast \in \{\text{full}, \text{restricted}\}$ to generate the RB spaces in each case by retaining  $50, 30, 35$ and $15$ bases functions for components $u_1, p_1, u_2$ and $p_2$, respectively. 
Then, in both cases, the ANN with 5 hidden layers with $16, 64, 128, 128$ and $128$ neurons, respectively, is constructed for each state subcomponent, using the hyperbolic tangent as an activation function. 
We use the same sets of snapshots used for the RB spaces construction to train the ANN and then validate the performance of the POD--NN method on the testing set of $26$ physical parameters as described above. Figures~\ref{fig:bfs_vacuum_mean_full} and~\ref{fig:bfs_vacuum_mean_restricted} show the relative errors of the POD--NN method in time, averaged over the training set of physical parameters, for each state subcomponent of the problem over the full and restricted time intervals, respectively.  It can be easily seen that in both cases the POD--NN is able to provide very good approximation for the velocity fields $u_1$ an $u_2$. On the other hand, the time discontinuity at the initial time step pervades in time and compromises the POD--NN approximation for the pressure fields $p_1$ and $p_2$, which does not happen if we restrict the analysis over the time interval where the discontinuity vanishes.  

Given that the POD--NN approach seems to be sensitive to the time discontinuity, as evidenced above, we perform another test case on the same domain where we have a smooth transition from the zero initial condition by considering the inlet boundary condition $u_{in}(x,y) =  \left(w(y,t),  \\0\right)^T$ on $\Gamma_{in}$, where the function $w$ is given by
\begin{equation*}
    w(y,t) =
    \begin{cases}
         \bar U \cdot \frac{1}{2} \left( 1 - \cos\left({\frac{2\pi t}{0.4}}\right)\right)\frac{4}{9} (y-2)(5-y), \ y \in [2,5], \ t \leq 0.4, \\
         \bar U \cdot\frac{4}{9} (y-2)(5-y), \ y \in [2,5], \ t > 0.4. \\
    \end{cases}
\end{equation*}
\begin{figure}
    \centering
    \begin{subfigure}[b]{0.49\textwidth}
        \includegraphics[width=\textwidth, trim = {26 0 24 0}, clip]{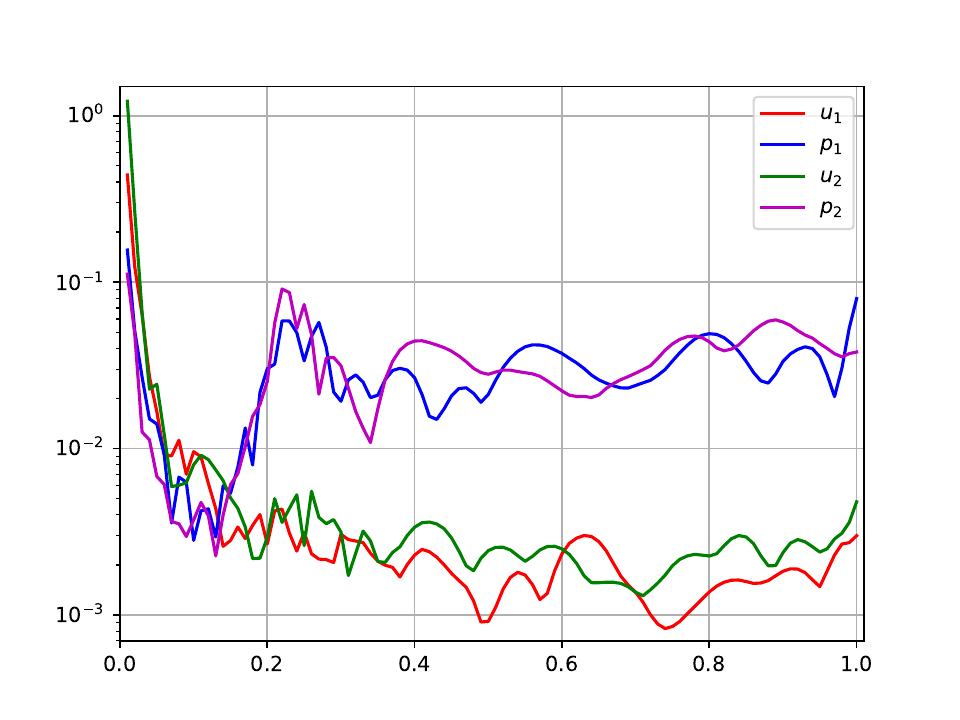}
        \caption{Mean relative errors over the entire time interval $[0,1]$}
         \label{fig:bfs_smooth_mean_full}
    \end{subfigure}
    \hfill
    \begin{subfigure}[b]{0.49\textwidth}
        \includegraphics[width=\textwidth, trim = {15 0 35 0}, clip]{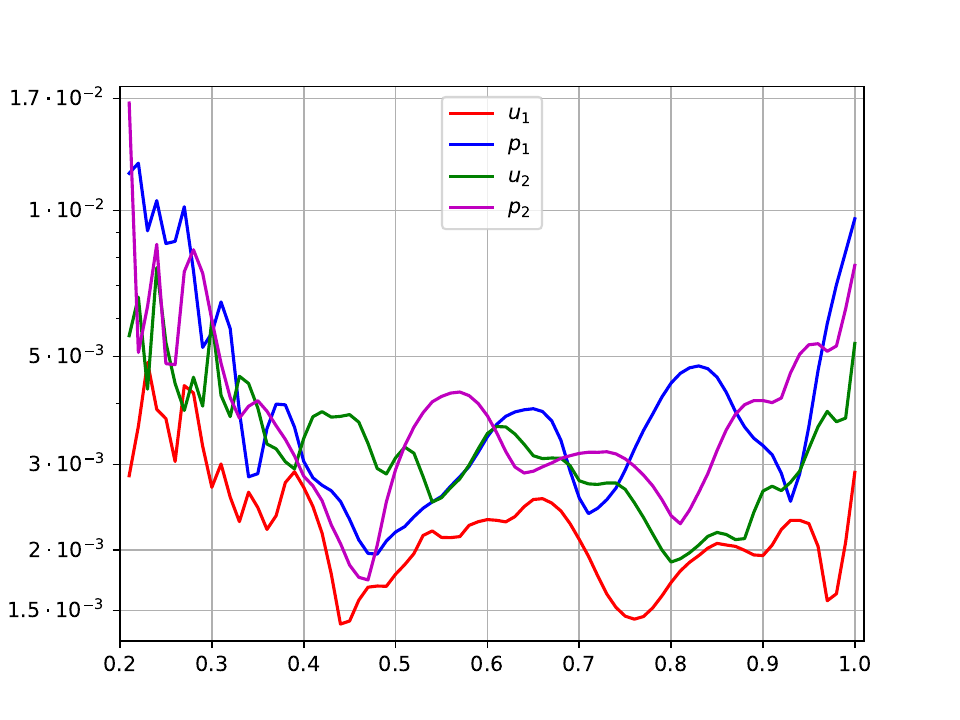}
         \caption{Mean relative errors over the restricted time interval $[0.21,1]$ }
         \label{fig:bfs_smooth_mean_restricted}
    \end{subfigure}    
    \caption{Mean over the parameter testing set of  POD--NN relative errors for the continuous in time backward--facing step test case  with respect to $\bar U$ over the entire (a) and restricted (b) time intervals }
    \label{fig:bfs_smooth_mean}
\end{figure}

As in the previous test case, we consider one physical parameter $\bar U$ to generate new POD--based reduced spaces and we train the ANN. Again, we consider 101 parameters, 75 in the training set and 26 in the test set, selected uniformly in the parameter space, see Table~\ref{table:offline_bfs}. 
The dimensions of the RB spaces and the architecture of the ANN are the same as in the discontinuous test case above. We again consider two scenarios by using use the sets $\mathcal{P}_{\text{train}}^{\text{full}}$ and $\mathcal{P}_{\text{train}}^{\text{restricted}}$ to perform POD compression and to train the ANN. Here, we discover another complication in the POD--NN approach: the sensitivity to small data, i.e., the ANN struggles to predict the data at the initial time steps as the norm is very close to zero. 
In Figures~\ref{fig:bfs_smooth_mean_full} and~\ref{fig:bfs_smooth_mean_restricted}, we show the relative errors of the POD--NN method in time, averaged over the training set of physical parameters, for each state subcomponent of the problem over the full and restricted time intervals, respectively. The picture here is very similar to the one in the discontinuous case, i.e. the velocity fields $u_1$ and $u_2$ show very good approximation properties (except for the first few time steps where the relative errors are expected to be bigger since the solution norm is very close to zero) while the POD--NN method struggles to learn the pressure field over the whole time intervals while the situation is improved by restricting the analysis to the time interval where the velocity fields norms are big enough. 

From the analysis carried out before, we can see that POD--NN approach provides excellent results as soon as we are far away enough from the ``problematic'' regions: either the regions that contain discontinuities or where the solution has a small norm. 
The error analysis reported in Figures~\ref{fig:bfs_vacuum_mean} and~\ref{fig:bfs_smooth_mean} shows that the relative error does not exceed the $1\%$--threshold if we exclude from the POD--NN model set--up the mentioned ``problematic'' regions, where the solution can be computed by FOM or POD--Galerkin methods with more reliability.
\begin{figure}
    \centering
    \includegraphics[width=0.65\textwidth]{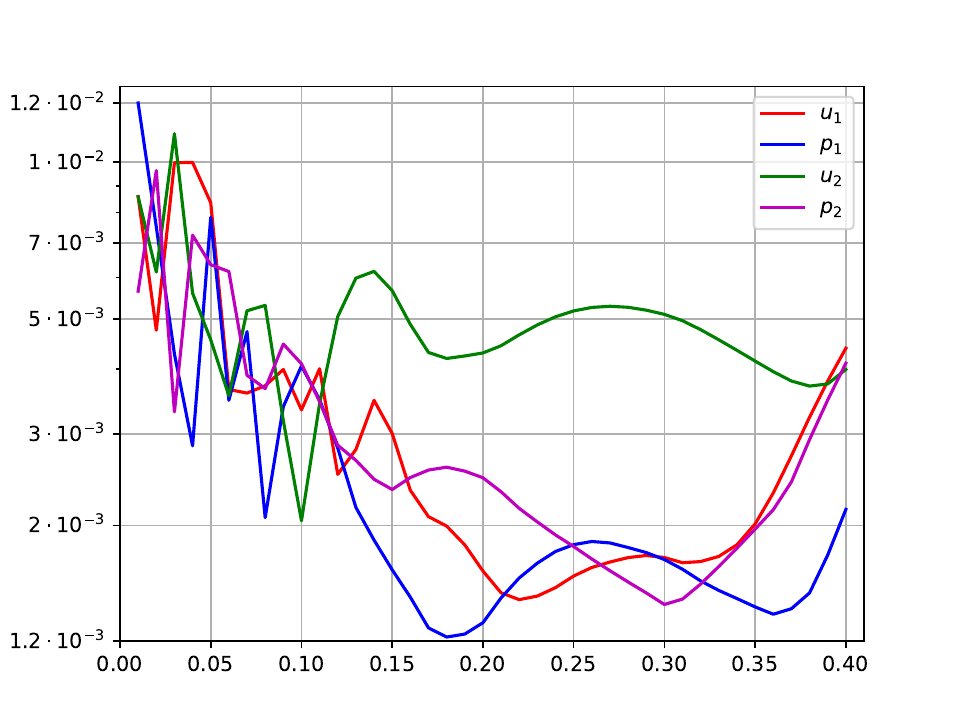}
    \caption{Mean over the parameter testing set of POD--NN relative errors for the discontinuous lid--driven cavity test case with respect to $\bar{U}$ over the time interval [0,0.4]  }
    \label{fig:cavity_vacuum_mean}
\end{figure}

Concerning the lid--driven cavity flow test case presented in~\ref{cavity_flow}, we were able to achieve a good approximation by NN--ROM method even in the presence of discontinuities, by merely extending the training set. As in the previous test cases, we draw $101$ uniformly sampled parameters from the parameter space reported in Table~\ref{table:offline_cavity},  $75$ for the training set and 26 for the testing set.
We perform the POD compression using the training set to generate RB spaces by retaining  $25, 15, 25$ and $15$ bases functions for components $u_1, p_1, u_2$ and $p_2$, respectively.
Then, the ANN with 5 hidden layers with $16, 64, 128, 128$ and $128$ neurons, respectively, is constructed for each state subcomponent, using the hyperbolic tangent as an activation function. 
We use the same $75$ snapshots used for the RB spaces construction to train the ANN and then validate the performance of the POD--NN method on the testing set of $26$ parameters as described above. 
In this case, the statistical analysis of the relative errors shows good behaviour on the entire time interval, as is evidenced in Figure~\ref{fig:cavity_vacuum_mean}.

The tests carried out in this section suggest the necessity to set up a proper sensitivity analysis with respect to model parameters and the problem geometry which will be the subject of future work.
}

\section{Conclusions}
\label{conclusions}

In this work, we described and conducted the convergence analysis of an optimisation--based domain decomposition algorithm for nonstationary Navier--Stokes equations.

The original problem cast into the optimisation--based domain--decomposition framework leads to the optimal control problem aimed at minimising the coupling error at the interface; the problem, then, has been tackled using an iterative gradient--based optimisation algorithm, which allowed us to obtain a complete separation of the solvers on different subdomains. 

At the reduced--order level, we provided two techniques: a POD--Galerkin projection and a data--driven POD--NN, both of them on separate domains. In the Galerkin projection, the optimal--control problem was solved much faster, not only because of the reduced dimensions but also because of the smaller number of iterations. In the POD--NN results are less accurate, but the computational time is way cheaper with respect to the other methods. \reviewerA{Additional statistical analysis has been carried out to investigate the approximation properties of POD--NN methodology which confirms the sensitivity of the method to discontinuities in time and small data. In future works, we want to perform a sensitivity analysis of the method to the model parameters and the problem geometry. }

As has been mentioned in the paper, the aforementioned techniques could be promising for various areas of computational physics. 
First of all, these algorithms can be used when complex time--dependent problems arise and domain decomposition becomes necessary due to the number of degrees of freedom.
\reviewerB{We plan to extend this algorithm also to several subdomains for more complex engineering problems as similarly done in \cite{taddei2023non}.}
Moreover, in the context of multi--physics or fluid--structure interaction, the coupling of pre--existing solvers on each subcomponent can be exploited in this framework, with the additional benefit of the reduction for parametric problems, which guarantees high adherence with respect to the full order solutions.
Finally, in case the codes are not directly accessible, the presented non--intrusive approach can be used to highly speed up the simulations while still obtaining meaningful results.

\section*{Acknowledgements}
This work was supported by the European Union's Horizon 2020 research and innovation programme under the Marie Sklodowska--Curie Actions [grant agreement 872442] (ARIA, Accurate Roms for Industrial Applications)  and by PRIN ``Numerical Analysis for Full and Reduced Order Methods for Partial Differential Equations'' (NA--FROM--PDEs) project. This project has received funding from the European High-Performance Computing Joint Undertaking (JU) under grant agreement No 955558. IP has been funded by a SISSA fellowship within the projects TRIM [project code: G95F21001070005] and iNEST. MN has been funded by the Austrian Science Fund (FWF) through project F 65 ``Taming Complexity in Partial Differential Systems" and project P 33477. DT has been funded by a SISSA Mathematical fellowship within Italian Excellence Departments initiative by Ministry of University and Research.

\newpage
\appendix
\section{Lifting Function and Velocity Supremiser Enrichment}
\label{lifting_supremiser}
The use of lifting functions is quite common in the reduced basis method (RBM) framework; see, for example, \cite{BallarinManzoniQuarteroniRozza2015, Rozza_book}. It is motivated by the fact that, in the chosen model, we tackle a problem with non--homogeneous Dirichlet boundary conditions on the parts of the boundaries $\Gamma_{i,D}, \,i=1,2$. From the implementation point of view, this does not present any problem when dealing with the high--fidelity model, since there are several well--known techniques for non--homogeneous essential conditions, in particular at the algebraic level. 
However, these boundary conditions create some problems when dealing with the reduced basis methods. 
Indeed, we seek to generate a linear vector space that is obtained by the compression of the set of snapshots, and this clearly cannot be achieved by using snapshots that satisfy different Dirichlet conditions, as the resulting space would not be linear. 
This problem is solved by introducing a lifting function $l_{i,h}^n \in V_{i,h},$ for $i =1,2,$ during the offline stage, such that $l_{i,h}^n = u_{i, D, h}^n$ on $\Gamma_{i,D}$. 
We define two new variables $u_{i,0,h}^n:=u_{i,h}^n-l_{i,h}^n \in V_{i,0,h },$ for $i=1,2,$ which satisfy the homogeneous condition $u_{i,0,h}^n = 0$ on $\Gamma_{i,D}$. 
So, they can be used to generate the reduced basis linear space. 
We remark that the lifting function is needed only in the domain where the Dirichlet boundary is non--empty, i.e., where $\Gamma_{i,D} \neq \emptyset$ for $i=1,2$.
It is important to point out that the choice of lifting functions is not unique; in our work, we chose to use the solution of the incompressible Stokes problem in one of the domains $\Omega$, $\Omega_1$ or $\Omega_2$ (depending on the particular model we are investigating) with the velocity equal to $u_{D}$ on the corresponding parts of the boundaries and the homogeneous Neumann conditions analogous to the original problem setting.   

Velocity supremiser enrichment is a very important technique to ensure that the reduced spaces are inf--sup stable in the context of saddle--point problems.  We recall that each velocity snapshot, which is a solution to the incompressible Navier--Stokes equation, is divergence--free. Hence, the term $b_i(\cdot, \cdot)$ for $i=1,2$ applied to any pair of functions in the span of the snapshots will be zero. This does not allow us to fulfil the inf--sup condition of the type \eqref{eq:infsup0}. For this reason, there is a need to enrich the reduced velocity spaces with extra functions, which are called supremisers, that will make the pairs of velocity--pressure reduced spaces inf--sup stable. 
The supremiser variables $s_{i,h}^n$, for $i=1,2,$ are defined as the solution to the following problem: find $s_{i,h}^n \in V_{i,0,h}$ such that 
\begin{equation}
    \left( \nabla v_{i,h},\nabla s_{i,h}^n\right) = b_{i,h} \left( v_{i,h}, p_{i,h}^n \right) \quad \forall v_{i,h} \in V_{i,0,h}, \label{eq:supremiser}
\end{equation}
where $p_{i,h}^n,$ for $i=1,2,$ are the finite--element pressure solutions of the Navier--Stokes problem and the left--hand side is the scalar product that defines the norm with which the variational spaces $V_{i,0,h}$ are endowed. For more details, we refer to \cite{BallarinManzoniQuarteroniRozza2015, GernerVeroy2012}.

\section{Reduced Basis Generation}
\label{POD}
Once we obtain the homogenised snapshots $u_{i,0,h}$ as described in Appendix~\ref{lifting_supremiser} and the \both{velocity supremisers } $s_{i,h}$ for $i=1,2$, we are ready to construct a set of reduced basis functions. A very common choice when dealing with Navier--Stokes equations is to use the POD technique; see, for instance, \cite{Rozza_book}. 
In order to implement this technique, we will need two main ingredients: the matrices of the inner products and the snapshot matrices, obtained by a full--order model (FOM) discretization as the one presented in the previous sections.
First, we define the basis functions for the FE element spaces used in the weak formulation \eqref{eq:functional_fem}, \eqref{eq:state_fem} and \eqref{eq:adjoint_fem}: $V_{i,0,h} = \text{span}\lbrace \phi_j^{u_i} \rbrace_{j=1}^{\mathcal N_h^{u_i}}$, $Q_{i,h} = \text{span}\lbrace \phi_j^{p_i} \rbrace_{j=1}^{\mathcal N_h^{p_i}}$ and $X_h=\text{span}\lbrace \phi_j^{g} \rbrace_{j=1}^{\mathcal N_h^{g}}$,
where $\mathcal N_h^{*},$ for $ * \in \left\{ u_1, p_1, u_2, p_2, g\right\}$, denotes the dimension of the corresponding FE space. 

We proceed by building the snapshot matrices. 
First, we sample the parameter space and draw a discrete set of $K$ parameter values.
Then, the snapshots, i.e.,  the high--fidelity FE solutions at each parameter value in the sampling set and at each time--step $t_1, ...,t_M $, are collected
into snapshot matrices $\mathcal S_{u_i} \in \mathbb{R}^{\mathcal N_h^{u_i}\times MK}$,  $\mathcal S_{s_i} \in \mathbb{R}^{\mathcal N_h^{u_i}\times MK}$, $\mathcal S_{p_i} \in \mathbb{R}^{\mathcal N_h^{p_i}\times MK}$,  for $i=1,2$ and $\mathcal S_{g} \in \mathbb{R}^{\mathcal N_h^g\times MK}$ for the corresponding values. 

The next step is to define the inner--product matrices $X_{u_i}$, $X_{p_i}$, for $i=1,2,$ and $X_g$:
\begin{align*}
 (X_{s_i})_{jk}=(X_{u_i})_{jk} & =  \left( \nabla \phi_{k}^{u_i}, \nabla \phi_{j}^{u_i} \right)_{\Omega_i}, \quad &\text{for } j, k = 1,..., \mathcal{N}_h^{u_i}, \ i=1,2, &\\
 (X_{p_i})_{jk} & =  \left(  \phi_{k}^{p_i},  \phi_{j}^{p_i} \right)_{\Omega_i}, \quad &\text{for } j, k = 1,..., \mathcal{N}_h^{p_i}, \ i=1,2, &\\
 (X_{g})_{jk} & =  \left(  \phi_{k}^{g},  \phi_{j}^{g} \right)_{\Gamma_0}, \quad &\text{for } j, k = 1,..., \mathcal{N}_h^{g}.&
\end{align*}

We are now ready to introduce the correlation matrices $\mathcal{C}_{u_i}$, $\mathcal{C}_{s_i}$, $\mathcal{C}_{p_i}$ for $i=1,2$ and $\mathcal{C}_g$, all of dimension $MK\times MK$, as:
\begin{eqnarray*}
\mathcal{C}_{*} := \mathcal S_{*}^T X_{*} S_{*}
\end{eqnarray*}
for every $* \in \{u_1, p_1, u_2, p_2, s_1, s_2, g\}$.

Once we have built the correlation matrices, we are able to carry out a POD compression on the sets of snapshots. This can be achieved by solving the following eigenvalue problems:
\begin{eqnarray}
\mathcal{C}_{*}\mathcal{Q}_{*} = \mathcal{Q}_{*} \Lambda_{*}  \label{eq:eigenvalue_problem}
\end{eqnarray}
where $* \in \{u_1, s_1, p_1, u_2, s_2, p_2, g\}$, $\mathcal Q_{*}$ is the eigenvectors matrix and $\Lambda_{*}$ is the diagonal eigenvalues matrix with eigenvalues ordered by decreasing order of their magnitude. The $k$--th reduced basis function for the component ${*}$ is then obtained by applying the matrix $\mathcal S_{*}$ to $\underline v_k^{*}$, the $k$--th column vector of the matrix $\mathcal{Q}_{*}$:
\begin{equation*}
    \Phi_k^{*}:=\frac{1}{\sqrt{\lambda_k^{*}}} \mathcal{S}_{*} \underline v_k^{*},  
\end{equation*}
where $\lambda_k^{*}$ is the $k$--th eigenvalue from \eqref{eq:eigenvalue_problem}. 
Therefore, we are able to form the set of reduced basis as
\begin{equation*}
\mathcal A^{*} := \left\{ \Phi_1^{*}, ..., \Phi_{N_\ast}^{*}   \right\},
\end{equation*}
where the integer numbers $N_{*}$ indicate the number of the basis functions used for each component for $\ast \in \{u_1, p_1, u_2, p_2, s_1, s_2, g\}$.
Now, it is time to include the supremiser enrichment of the velocities spaces discussed at the beginning of this section. We provide the following renumbering of the functions for further simplicity: 
\begin{equation*}
    \Phi_{N_{u_i}+k}^{u_i}:=\Phi_{k}^{s_i}, \quad \text{for} \ k=1,...,N_{s_i}, \ i=1,2,
\end{equation*}
and we redefine $N_{u_i}:=N_{u_i}+N_{s_i}$, and new basis functions sets
\begin{equation*}
\mathcal A^{u_i} := \left\{ \Phi_1^{u_i}, ..., \Phi_{N_{u_i}}^{u_i}   \right\},
\end{equation*}
for $i=1,2$ and these new sets are now including extra basis functions obtained from the corresponding supremiser.
Finally, we introduce three separate reduced basis spaces -- for the state and the control variables, respectively:
\begin{eqnarray*}
    V_{rb}^{\ast} = \text{span} (\mathcal{A}^{\ast}), & & \text{dim}(V_{rb}^*) = N_{\ast},
\end{eqnarray*}
for $\ast \in \{u_1, p_1, u_2, p_2, g\}$. Now, due to the supremiser enrichment the spaces $V_{rb}^{u_i}$ and $V_{rb}^{p_i}$ are inf--sup stable in the sense \eqref{eq:infsup0} for $i=1,2$; the proof can be found in \cite{BallarinManzoniQuarteroniRozza2015}.
\newpage
\bibliographystyle{abbrv} 
\bibliography{main}

\end{document}